\newif\ifecma
\ecmafalse

\ifecma
\documentclass[12pt]{article}
\else
\documentclass{article}
\fi

\usepackage{amsmath, amsthm, amssymb}
\usepackage{graphicx}
\usepackage{verbatim}
\usepackage{natbib}
\usepackage{caption}
\usepackage{subcaption}
\usepackage{fancyvrb}
\usepackage{enumerate}
\usepackage{relsize}
\usepackage{tikz}
\usepackage{tikz-qtree}
\usepackage{setspace}

\usepackage{multirow}

\usepackage{titling}

\usepackage{hyperref}

\ifecma
\usepackage[margin=1in]{geometry}
\else
\usepackage[margin=1.46in]{geometry}
\fi

\hypersetup{colorlinks,citecolor=blue,urlcolor=blue,linkcolor=blue}

\usepackage[colorinlistoftodos]{todonotes}

\usepackage{stefan_tex}
\graphicspath{{./figures/}}

\theoremstyle{plain}
\newtheorem{prop}{Proposition}

\newtheorem{coro}[prop]{Corollary}
\newtheorem{lemm}[prop]{Lemma}
\newtheorem{theo}[prop]{Theorem}

\theoremstyle{definition}

\newtheorem{assu}{Assumption}

\theoremstyle{remark}

\newcommand\blfootnote[1]{%
  \begingroup
  \renewcommand\thefootnote{}\footnote{#1}%
  \addtocounter{footnote}{-1}%
  \endgroup
}

\author{Susan Athey \\ \texttt{athey@stanford.edu}
\and
Stefan Wager \\ \texttt{swager@stanford.edu}}
\date{\ifcase\month\or
January\or February\or March\or April\or May\or June\or
July\or August\or September\or October\or November\or December\fi \ \number%
\year\ \  }
\title{Policy Learning with Observational Data}

\begin{document}

\maketitle

\begin{abstract}
In many areas, practitioners seek to use observational data to learn a treatment assignment policy that satisfies application-specific constraints, such as budget, fairness, simplicity, or other functional form constraints.  For example, policies may be restricted to take the form of decision trees based on a limited set of easily observable individual characteristics. We propose a new approach to this problem motivated by the theory of semiparametrically efficient estimation. Our method can be used to optimize either binary treatments or infinitesimal nudges to continuous treatments, and can leverage observational data where causal effects are identified using a variety of strategies, including selection on observables and instrumental variables. Given a doubly robust estimator of the causal effect of assigning everyone to treatment, we develop an algorithm for choosing whom to treat, and establish strong guarantees for the asymptotic utilitarian regret of the resulting policy. 

\vspace{\baselineskip}

\end{abstract}

\ifecma
\setstretch{1.4}
\fi

\section{Introduction}

The problem of learning treatment assignment policies, or mappings from individual characteristics
to treatment assignments, is ubiquitous in applied economics and
statistics.\blfootnote{Forthcoming in Econometrica.
We circulated an earlier draft of this paper under the title
``Efficient Policy Learning''; the current title was chosen following guidance from the review process.
We are grateful for helpful conversations with colleagues
including
Victor Chernozhukov,
David Hirshberg,
Guido Imbens,
Michael Kosorok,
Alexander Luedtke,
Eric Mbakop,
Whitney Newey,
Xinkun Nie,
Molly Offer-Westort,
Alexander Rakhlin,
James Robins,
Erik Sverdrup,
Max Tabord-Meehan and
Zhengyuan Zhou,
and for feedback from the editor, referees, as well as seminar participants at a variety of universities and workshops.
We thank Guido Imbens for sharing the GAIN dataset with us. Financial support was provided by the Sloan Foundation,
Office of Naval Research grant N00014-17-1-2131,
National Science Foundation grant DMS-1916163,
and a Facebook Faculty Award.} It arises,
for example, in medicine when a doctor must decide which patients to refer
for a risky surgery; in marketing when a company needs to
choose which customers to send targeted offers to; and in government and policy settings,
when assigning students to educational programs or inspectors to buildings and restaurants.

The treatment assignment problem rarely arises in an unconstrained environment. Treatments are often expensive,
and so a policy may need to respect budget constraints.  Policies may need to be implemented in environments characterized
by human or machine constraints; for example, emergency medical professionals or police officers may need to implement decision
policies in the field, where a simple decision tree might be used. For internet or mobile services, algorithms may need to determine the set of information displayed to a user very quickly, and a simple lookup table may decrease the time it takes to respond to a user's request.
Fairness constraints may require a treatment assignment policy to depend only on particular types of covariates (for example, test scores or income),
even when other covariates are observed.

This paper is about using observational data to learn policies that respect the types of constraints outlined above.
The existing literature on policy learning has mostly focused on the setting where we want
to optimize allocation of a binary treatment using data from a randomized trial, or from a study with a known, random treatment assignment policy. 
In many problems, however, one may need to leverage richer forms of observational
data to learn treatment assignment rules. For example, if we want to learn whom to prescribe a drug
to based on data from a clinical trial, we need to have methods that deal with non-compliance
and resulting endogenous treatment assignments.\footnote{If we believed that compliance
patterns when we deploy our policy would be similar to those in the clinical trial, then an
intent-to-treat analysis may be a reasonable way to side-step endogeneity concerns. However,
if we suspect that compliance patterns may change (e.g., if patients may be more likely to adhere
to a treatment regime prescribed by their doctor than one randomly assigned in a clinical trial),
then using an analysis that disambiguates received treatment from assigned treatment is
necessary.} Or, if we are interested in offering some customers discounts, then we need methods
that let us study interventions to continuous variables (e.g., price) rather than just discrete ones.
The goal of this paper is to develop methods for policy learning that don't just work in randomized
trials (or related settings), but can instead work with a rich variety of observational designs.

Formally, we study the problem where we have access to observational data
and want to use it to learn a policy that maps a subject's characteristics
$X_i \in \xx$ to a binary decision, $\pi : \xx \rightarrow \cb{0, \, 1}$. The practitioner has also specified
a class $\Pi$ that encodes problem-specific constraints pertaining to budget, functional form, fairness, etc., and
requires that our learned policy \smash{$\hpi$} satisfies these constraints, \smash{$\hpi \in \Pi$}.
Then, following \citet{manski2004statistical,manski2009identification},
\citet{hirano2009asymptotics}, \citet{stoye2009minimax,stoye2012minimax} and \citet{kitagawa2015should},
we seek guarantees on the regret \smash{$R(\hpi)$}, i.e., the difference between the expected utility from deploying
the learned policy \smash{$\hpi$} over a target population and the best utility that could be achieved from deploying
any policy in the class $\Pi$ over the population.

Our paper builds on a rich literature at the intersection of econometrics,
statistics and computer science on learning structured treatment assignment rules, including
\citet{kitagawa2015should}, \citet{swaminathan2015batch} and \citet*{zhao2012estimating}.
Most closely related to us, \citet{kitagawa2015should} study a special case
of our problem where treatments are binary and exogenous with known assignment probabilities, and show that an algorithm
based on inverse-probability weighting achieves regret that depends optimally on the sample size and the
complexity of the policy class $\Pi$.\footnote{\citet{kitagawa2015should} also consider the case where
treatment assignment probabilities are unknown; in this case, however, their method no longer achieves
optimal dependence on the sample size.}

Here, we develop a new family of algorithms that achieve regret guarantees
with optimal dependence on sample size and on $\Pi$, but under considerably more generality on the sampling design.
We consider both the classical case where we want to optimize a binary treatment, and a related setting
where we want to optimize infinitesimal nudges to a continuous treatment (e.g., a price). Moreover,
our approach can leverage observational data where the treatment assignment mechanism may either be
exogenous with unknown assignment probabilities, or endogenous, in which case we require an instrument. 

Our approach starts from recent unifying results of \citet*{chernozhukov2016locally} on semiparametrically
efficient estimation. As discussed in more detail in Section \ref{sec:identification}, \citet{chernozhukov2016locally} show
that in many problems of interest, we can construct efficient estimates of average-treatment-effect-like
parameters $\theta$ as
\begin{equation}
\label{eq:dr_tau}
\htheta = \frac{1}{n} \sum_{i = 1}^n \hGamma_i,
\end{equation}
where \smash{$\hGamma_i$} is an appropriate doubly robust score
for the target estimand under the intervention of interest. This approach
can be used to target the average effect of a binary treatment, the average derivative of a continuous treatment,
and other related estimands.

In this paper we find that, whenever one can estimate the average utility of treating everyone\footnote{Throughout
this paper, we assume that there is no interference, i.e., assigning one unit to treatment doesn't affect outcomes
for others. For a discussion of treatment effect estimation under intereference, see \citet{hudgens2008toward},
\citet{manski2013identification}, and references therein.} using an estimator
of the type \eqref{eq:dr_tau} built via the doubly robust construction of \citet{chernozhukov2016locally}, we can
also usefully learn whom to target with the intervention via a simple procedure:
Given a pre-specified policy class $\Pi$ (e.g., linear decision rules or finite-depth decision trees),
we propose using the treatment assignment rule $\hpi$ that solves\footnote{If
this optimization problem has multiple solutions, we set $\hpi$ to an arbitrary maximizer of the objective. Our
formal results apply simultaneously to all solutions of \eqref{eq:meta_estimator}.}
\begin{equation}
\label{eq:meta_estimator}
\hpi = \argmax\cb{\frac{1}{n} \sum_{i = 1}^n \p{2\pi(X_i) - 1} \hGamma_i : \pi \in \Pi},
\end{equation}
where \smash{$\hGamma_i$} are the same doubly robust scores as used in \eqref{eq:dr_tau}. 
Our main result is that, under regularity conditions, the resulting policies \smash{$\hpi$} have
regret \smash{$R(\hpi)$} bounded on the order of \smash{$\sqrt{\VC(\Pi)/n}$} with high probability.
Here, $\VC(\Pi)$ is the Vapnik-Chervonenkis dimension of the
class $\Pi$ and $n$ is the sample size. We also highlight how the constants in this bound depend
on fundamental quantities from the semiparametric efficiency literature.

Our proof combines results from semiparametrics with carefully tailored analysis tools that
build on classical ideas from empirical process theory. The reason we obtain strong guarantees for the
approach \eqref{eq:meta_estimator} is closely tied to robustness properties of the estimator \eqref{eq:dr_tau}.
In the setting where we only want to estimate a single average effect parameter, it is well known that non-doubly
robust estimators can also be semiparametrically efficient \citep*{hirano2003efficient}. Here, however, we
need convergence results that are strong enough to withstand optimization over the whole class $\Pi$.
The fact that doubly robust estimators are fit for this task is closely related to their ability to achieve semiparametric
efficiency under general conditions, even if nuisance components are estimated via black-box machine learning
methods for which we can only guarantee fast enough convergence in
mean-squared error \citep{chernozhukov2016double,van2011targeted}.

We spell out our general framework in Section \ref{sec:method}. For intuition, however, it is helpful
to first consider this approach in the simpler case where we want to study the effect of a binary
treatment $W_i \in \cb{0, \, 1}$ on an outcome $Y_i \in \RR$ interpreted as a utility and are willing to assume selection on observables (unconfoundedness):
We have potential outcomes $\cb{Y_i(0), \, Y_i(1)}$ such that $Y_i = Y_i(W_i)$ and $\cb{Y_i(0), \, Y_i(1)} \indep W_i \cond X_i$
\citep{imbens2015causal}. Then, the utilitarian regret of deploying a policy \smash{$\pi \in \Pi$} is \citep{manski2009identification}
\begin{equation}
\label{eq:binary_regret}
R(\pi) = \max\cb{\EE{Y_i(\pi'(X_i))} : \pi' \in \Pi} - \EE{Y_i(\pi(X_i))},
\end{equation}
and we can construct our estimator \eqref{eq:meta_estimator} using
the well known augmented inverse-propensity weighted scores of \citet*{robins1994estimation},\footnote{See
Section \ref{sec:gain} for a detailed discussion of how to implement our policy learner \eqref{eq:meta_estimator}
based on these augmented inverse-propensity weighted scores in practice.}
\begin{equation}
\label{eq:aipw_score}
\begin{split}
&\hGamma_i = \hatm(X_i, \, 1) - \hatm(X_i, \, 0) + \frac{W_i - \he(X_i)}{\he(X_i)(1 - \he(X_i))}\p{Y_i - \hatm(X_i, \, W_i)}, \\
&e(x) = \PP{W_i  = 1\cond X_i = x}, \ \ m(x, \, w) = \EE{Y_i(w) \cond X_i = x},
\end{split}
\end{equation}
where \smash{$\he(x)$} and \smash{$\hatm(x, \, w)$} denote non-parametric estimates of $e(x)$ and $m(x, \, w)$ respectively.
In this setup, our result implies that---under regularity conditions---the estimator \eqref{eq:meta_estimator} with
scores \eqref{eq:aipw_score} has regret \eqref{eq:binary_regret} bounded on the order of \smash{$\sqrt{\VC(\Pi)/n}$}.

Even in this simplest case, our result is considerably stronger than results currently available in the
literature. The main result of \citet{kitagawa2015should} is that, if treatment propensities $e(X_i)$ are known, then
a variant of inverse-propensity weighted policy learning achieves regret on the order of \smash{$\sqrt{\VC(\Pi)/n}$}.
However, in observational studies where the treatment propensities are unknown, the bounds of
\citet{kitagawa2015should} depend on the rate at which we can estimate $e(\cdot)$, and will generally
decay slower than $1/\sqrt{n}$. The only other available $1/\sqrt{n}$-bounds for policy learning in
observational studies with a binary treatment that we are aware of are a result of \citet*{van2006cross}
for the case where $\Pi$ consists of a finite set of policies whose cardinality grows with $n$, and a
result of \citet{kallus2017balanced} in the special case $m(\cdot, \, w)$ is assumed to belong to a
reproducing kernel Hilbert space. The idea of using doubly robust scores to learn optimal treatment
assignment of a binary treatment has been previously discussed in \citet*{langford2011doubly} and
\citet*{zhang2012estimating}; however, neither paper provides a regret bound for this approach.

In the more general case where the observed treatment assignments $W_i$ may be continuous and/or
we may need to use instrumental variables to identify causal effects, both the methods and regret
bounds provided here are new.  By connecting the policy learning problem to the semiparametric efficiency literature,
we are able to develop a general framework that applies across a variety of settings.

\subsection{Related Work}

The literature on optimal treatment allocation has been
rapidly expanding across several fields. In the econometrics literature, the program of learning regret-optimal treatment
rules was started by \citet{manski2004statistical,manski2009identification}. One line of work considers the case where the policy class
is unrestricted, and the optimal treatment assignment rule simply depends on the sign of the conditional average treatment effect for each individual unit. In this setting,
\citet{hirano2009asymptotics} show that when $1/\sqrt{n}$-rate estimation of the conditional
average treatment effect function is possible, then treatment assignment rules obtained by thresholding
an efficient estimate of the conditional average treatment effect are asymptotically minimax-optimal.
Meanwhile, \citet{stoye2009minimax} derives finite sample minimax decision
rules in a class of problems where both the response surfaces and the policies $\pi$ may depend arbitrarily on
covariates. Further results are given in \citet{armstrong2015inference}, \citet{bhattacharya2012inferring},
\citet{chamberlain2011bayesian}, \citet{dehejia2005program}, \citet{kasy2016partial}, \citet{stoye2012minimax} and
\citet{tetenov2012statistical}.

Building on this line of work, \citet{kitagawa2015should} study policy learning in a non-parametric setting where the learned
policy $\hpi$ is constrained to belong to a structured class $\Pi$ and show that, in this case, we can obtain regret bounds
relative to the best policy in $\Pi$ that scale with the complexity of the class $\Pi$.
A key insight from \citet{kitagawa2015should} is that, when propensity scores are known and $\Pi$ has
finite VC dimension, it is possible to get $1/\sqrt{n}$-rate regret bounds for policy learning over a class
$\Pi$ even if the conditional average treatment effect function itself cannot be estimated at a $1/\sqrt{n}$-rate;
in other words, we can reliably find a nearly best-in-class policy without needing to accurately estimate a model that
describes all causal effects. As discussed above, our paper builds on this work by considering rate-optimal
regret bounds for best-in-class policy learning in observational studies where propensity scores are unknown
and treatment assignment may be endogenous, etc.

One difference between our results and those of \citet{kitagawa2015should} is that the latter provide
finite sample regret bounds, whereas our results are asymptotic in the sample size $n$. The reason
for this is that our bounds rely on results from the literature on semiparametric estimation
\citep*{bickel,chernozhukov2016locally,chen2008semiparametric,hahn1998role,newey1994asymptotic,robins1},
which themselves are asymptotic. Recently, \citet{armstrong2017finite} showed that, in a class of average
treatment effect estimation problems, finite sample conditionally minimax linear estimators are asymptotically
efficient, thus providing a connection between desirable finite sample guarantees and asymptotic optimality.
It would be interesting to examine whether similar connections are possible in the policy learning case.

Policy learning from observational data has also been considered in parallel literatures developed in both statistics
\citep*{luedtke2016statistical,qian2011performance,zhang2012estimating,zhao2012estimating}
and machine learning
\citep*{beygelzimer2009offset,langford2011doubly,kallus2017balanced,swaminathan2015batch}.
Two driving themes behind these literatures are the development of performant algorithms for solving
the empirical maximization problems (and relaxations thereof) that underlie policy learning, and the use of
doubly robust objectives for improved practical performance. \citet{kallus2017balanced}, \citet{swaminathan2015batch} and
\citet{zhao2012estimating} also prove regret bounds for their methods; however,
they do not achieve a $1/\sqrt{n}$ sample dependence, with the exception of \citet{kallus2017balanced} in the special case of the reproducing kernel
Hilbert space setting described above. Finally, \citet{luedtke2017faster} propose a class of regret bounds
that decay faster than $1/\sqrt{n}$ by exploiting non-uniform asymptotics; see Section \ref{sec:lb} for a further discussion.

The problem of optimal treatment allocation can also be seen as a special case of the broader problem
of optimal data-driven decision making. From this perspective, our result is related to the work of
\citet{rudin2018big} and \citet{bertsimas2014predictive}, who study data-driven rules for optimal inventory management and related problems.
Much like in our case, they advocate learning with a loss function that is directly tied to a utility-based criterion.
Finally, we note a growing literature on estimating conditional average treatment effects, including
\citet{athey2015machine}, \citet*{athey2018generalized}, \citet{nie2017learning}, and references therein.
Although the goal is similar to that of learning optimal treatment assignment rules, the
specific results themselves differ; they focus on squared-error loss rather than utilitarian regret.

\section{From Efficient Policy Evaluation to Learning}
\label{sec:method}

Our goal is to learn a policy $\pi \in \Pi$ that maps a subject's features $X_i \in \xx$
to a treatment decision: $\pi : \xx \rightarrow \cb{0, \, 1}$. In order to do so, we assume
that we have independent and identically distributed samples $(X_i, \, Y_i, \, W_i, \, Z_i)$,
where $Y_i \in \RR$ is the outcome we want to intervene on,
$W_i$ is the observed treatment assignment, and $Z_i$ is an (optional) instrument used for identifying causal effects.
In cases where $W_i$ is exogenous, we simply take $Z_i = W_i$.
Throughout our analysis, we interpret $Y_i$ as the utility resulting from our intervention on
the $i$-th sample,  e.g., $Y_i$ could measure the benefit accrued by a subject minus
a potentially personalized cost of treatment (in Section \ref{sec:gain} we demonstrate inclusion
of linear costs in the context of an application). We then seek policies that make the expected value
of $Y_i$ large.
 
We define the causal effect of the intervention $\pi(\cdot)$
in terms of the potential outcomes model \citep{neyman1923applications,rubin1974estimating},
whereby the $\cb{Y_i(w)}$ correspond to utilities we would have observed for the $i$-th sample
had the treatment been set to $W_i = w$, and $Y_i = Y_i(W_i)$. When instruments are present,
we always assume that the exclusion restriction holds so that this notation is well specified.
We consider both examples with a binary treatment $W_i \in \cb{0, \,1}$ and with a continuous
treatment $W_i \in \RR$.

In the case where $W_i$ is binary, we follow the existing literature
\citep{hirano2009asymptotics,kitagawa2015should,manski2004statistical,stoye2009minimax},
and study interventions that directly specify the treatment level. In this case, the
utility of deploying a policy $\pi(\cdot)$ relative to treating no one is \citep{manski2009identification}
\begin{equation}
\label{eq:V_binary}
V(\pi) = \EE{Y_i(\pi(X_i)) - Y_i(0)},
\end{equation}
and the corresponding policy regret relative to the best possible policy in the class $\Pi$ is
\begin{equation}
\label{eq:regret}
R(\pi) = \max\cb{V(\pi') : \pi' \in \Pi} - V(\pi).
\end{equation}
As discussed in the introduction, in this binary setting, \citet{kitagawa2015should} show that
if $W_i$ is exogenous with known treatment propensities, then we can use inverse-propensity
weighting to derive a policy \smash{$\hpi$} whose regret \smash{$R(\hpi)$} decays
as $1/\sqrt{n}$, with
\begin{equation}
\label{eq:IPW}
\hpi_{IPW} = \argmax\cb{\frac{1}{n} \sum_{i = 1}^n \frac{1\p{\cb{W_i = \pi(X_i)}} Y_i}{\PP{W_i = \pi(X_i) \cond X_i}} : \pi \in \Pi}.
\end{equation}
Here, we develop methods that can also be used
in observational studies where treatment propensities may be unknown, and where we
may need to use instrumental variables to identify $V(\pi)$ from \eqref{eq:V_binary}.

Meanwhile, when $W_i$ is continuous, we study infinitesimal interventions on the treatment level motivated by
the work of \citet*{powell1989semiparametric}. We define the utility of such an infinitesimal
intervention as
\begin{equation}
\label{eq:V_continuous}
V(\pi) = \sqb{\frac{d}{d\nu} \EE{Y_i(W_i + \nu \pi(X_i))}}_{\nu = 0},
\end{equation}
and then define regret in terms of $V(\pi)$ as in \eqref{eq:regret}. 
One interesting conceptual difference that arises in this case is that, now, our interventions $\pi(X_i) \in \cb{0, \, 1}$
and observed treatment assignments $W_i \in \RR$ may take values in different spaces. This can arise,
for example, if we want to target customers with personalized discounts
and have access to past prices $W_i$ that take on a continuum of values, but
are restricted to considering a class of interventions that only allow us to make a binary decision $\pi(X_i) \in \cb{0, \, 1}$
on whether to offer each customer a small discount or not. The fact that we can still learn low-regret
policies via the simple strategy \eqref{eq:meta_estimator} even when these two spaces are decoupled
highlights the richness of the policy learning problem.\footnote{Another interesting question one could
ask is how best to optimize the assignment of $W_i$ globally rather than locally (i.e., the case
where we can set the treatment level $w$ to an arbitrary level, rather than simply nudge the pre-existing
levels of $W_i$). This question would require different formal tools, however, as the results developed
in this paper only apply to binary decisions.}

With both binary and continuous treatments, the regret
of a policy $\pi$ can be written in terms of a conditional average treatment effect function,
\begin{equation}
\label{eq:cate}
\tau(x) = \EE{Y_i(1) - Y_i(0) \cond X_i = x} \ \text{ or } \ \tau(x) = \sqb{\frac{d}{d\nu} \EE{Y_i(W_i + \nu) \cond X_i = x}}_{\nu = 0},
\end{equation}
such that $V(\pi) = \EE{\pi(X_i) \tau(X_i)}$ and regret $R(\pi)$ is as in \eqref{eq:regret}.
Our analysis pertains to any setup with a regret function $R(\pi)$ that admits such a representation.
Given these preliminaries, recall that our goal is to learn low regret policies, i.e., to use observational data to derive a
policy $\hpi \in \Pi$ with a guarantee that $R(\hpi) = \oo_P\p{1/\sqrt{n}}$. In order to
do so, we need to make assumptions on the observational data generation distribution
that allow for identification and adequate estimation of $V(\pi)$, and also control the size
of $\Pi$ in a way that makes emulating the best-in-class policy a realistic objective.
The following two subsections outline these required conditions; our main result is
then stated in Section \ref{sec:main_result}.

\subsection{Identifying and Estimating Causal Effects}
\label{sec:identification}

In order to learn a good policy $\hpi$, we first need to be able to evaluate $V(\pi)$ for any specific policy $\pi$. 
Our main assumption, following \citet*{chernozhukov2016locally}, is that we can construct a
doubly robust score for the average treatment effect $\theta = \EE{\tau(X_i)}$.
At the end of this section we discuss how this approach applies to three important examples,
and refer the reader to \citet{chernozhukov2016locally} for a more general discussion of
when such doubly robust scores exist.

\begin{assu}
\label{assu:identification}
Write $m(x, \, w) = \EE{Y_i(w) \cond X_i = x} \in \mm$ for the counterfactual response surface.
We assume that $m(x, \, w)$ induces a treatment effect function $\tau_m(x, \, w)$ with the
following properties:
\begin{enumerate}
\item The functional $m(\cdot) \rightarrow \tau_m(\cdot)$ is linear in $m$, and there exists
a weighting function $g(x, \ z)$ that identifies $\tau_m(\cdot)$ via
\begin{equation}
\label{eq:representer}
\EE{\tau_{\widetilde{m}}(X_i, \, W_i) - g(X_i, \, Z_i)\widetilde{m}(X_i, \, W_i) \cond X_i} = 0,
\end{equation}
for any counterfactual response surface $\widetilde{m}(x, \, w) \in \mm$.
\item Policy value can be defined in terms of moments of $\tau_m(X_i, \, W_i)$, such that
 $V(\pi) = \EE{\pi(X_i)\tau(X_i)}$ with $\tau(x) = \EE{\tau_m(X_i, \, W_i) \cond X_i = x}$ for
 all $\pi : \xx \rightarrow \cb{0, \, 1}$.
\end{enumerate}
In some examples $\tau_m(x, \, w)$ does not depend on $w$, and we
omit the $w$-argument of $\tau_m(\cdot)$. 
\end{assu}

Given this setup, \citet{chernozhukov2016locally} propose first estimating $g(\cdot)$ and $m(\cdot)$,
and then consider
\begin{equation}
\label{eq:LR}
\htheta = \frac{1}{n} \sum_{i = 1}^n \hGamma_i, \ \ \hGamma_i = \tau_{\hatm}(X_i, \, W_i) + \hg\p{X_i, \, Z_i} \p{Y_i - \hatm\p{X_i, \, W_i}}.
\end{equation}
They show that this estimator is $\sqrt{n}$-consistent and asymptotically unbiased Gaussian for $\theta$,
provided that the nuisance estimates $\hg(\cdot)$ and $\hatm(\cdot)$
converge sufficiently fast and that we use cross-fitting \citep{chernozhukov2016double,schick1986asymptotically}.
This estimator is also semiparametrically efficient under general conditions
\citep{newey1994asymptotic}.\footnote{Our results don't depend on efficiency of \eqref{eq:LR}; rather,
we only use $\sqrt{n}$-consistency. In cases where \eqref{eq:LR} may not be efficient, our regret bounds
still hold verbatim; the only difference being that we can no longer interpret the terms of the form
\smash{$\EE{\Gamma_i^2}$} appearing in the bound as related to the semiparametric efficient variance for $\theta$.}

Our approach to policy learning builds on these foundations. We again start by estimating nuisance components
and by forming doubly robust scores as in \eqref{eq:LR}. However, instead of just averaging the \smash{$\hGamma_i$}
to estimate $\theta$, we use these scores for policy learning by plugging them into \eqref{eq:meta_estimator}.
Our main result will establish that we can get strong regret bounds for learning policies under conditions
that are similar to those used by \citet{chernozhukov2016locally} to show asymptotic normality of \eqref{eq:LR} and, more
broadly, that build on assumptions often made in the literature on semiparametric efficiency
\citep*{bickel,chen2008semiparametric,hahn1998role,newey1994asymptotic,robins1}.

As in the recent work of
\citet{chernozhukov2016double} on double machine learning or that of \citet{van2011targeted} on targeted learning,
we take an agnostic view on how the nuisance estimates $\hg(\cdot)$ and $\hatm(\cdot)$ are obtained,
and simply impose high level conditions on their rates of convergence.
Given sufficient regularity, we can construct estimators that satisfy the rate condition \eqref{eq:an_rate} via, e.g.,
sieve-based methods \citep*{chen2007large}
or kernel regression \citep*{caponnetto2007optimal}.
Moreover, in applications, we may want to consider several different machine learning methods for each component,
or potentially combinations thereof, and then use cross-validation to choose which method to use.
For completeness, we allow problem specific quantities to change with the sample size $n$,
and track this dependence with a subscript $n$, e.g., $m_n(x, \, w) = \EE[n]{Y_i(w) \cond X_i = x}$, etc.

\begin{assu}
\label{assu:LR}
In the setting of Assumption \ref{assu:identification}, assume that second moments are controlled as
\smash{$\EE[n]{m_n^2(X_i, \, W_i)}$}, \smash{ $\EE[n]{\tau_{m_n}^2(X_i, \, W_i)} < \infty$}
and \smash{$\EE[n]{g_n^2(X_i, \, Z_i)} < \infty$} for all $n = 1, \,2, \, ...$,
and that we have access to uniformly consistent estimators
of these nuisance components,
\begin{equation}
\begin{split}
&\sup_{x, \, w}\cb{\abs{\hatm_n(x, \, w) - m_n(x, \, w)}}, \
\sup_{x, \, w}\cb{\abs{\tau_{\hatm_n}(x, \ w) - \tau_{m_n}(x, \, w)}} \rightarrow_p 0, \\
&\sup_{x,\, z}\cb{\abs{\hg_n(x, \, z) - g_n(x, \, z)}} \rightarrow_p 0,
\end{split}
\end{equation}
whose $L_2$ errors decay as follows, for some $0 < \zeta_m, \, \zeta_g < 1$ with $\zeta_m + \zeta_g \geq 1$
and some $a(n) \rightarrow 0$, where $(X, \, W, \, Z)$ is taken to be an independent test
example drawn from the same distribution as the training data:\footnote{A notable special case of
this assumption is when $\zeta_m = \zeta_g = 1/2$; this is equivalent to the standard assumption
in the semiparametric estimation literature that all nuisance components (i.e., in our case, both the
outcome and weighting regressions) are $o(n^{-1/4})$-consistent in terms of $L_2$-error. The
weaker requirement \eqref{eq:an_rate} reflects the fact that doubly robust treatment effect
estimators can trade-off accuracy of the $m$-model with accuracy of the $g$-model, provided
the product of the error rates is controlled \citep{farrell2015robust}.}
\begin{equation}
\label{eq:an_rate}
\begin{split}
&\EE{\p{\hatm_n(X, \, W) - m_n(X, \, W)}^2}, \, \EE{\p{\tau_{\hatm_n}(X, \, W) - \tau_{m_n}(X, \, W)}^2} \leq \frac{a(n)}{n^{\zeta_m}}, \\
&\EE{\p{\hg_n(X, \, Z) - g_n(X, \, W)}^2} \leq \frac{a(n)}{n^{\zeta_g}}.
\end{split}
\end{equation}
\end{assu}

We end this section by verifying that Assumption \ref{assu:identification} in fact covers several settings of
interest, and is closely related to several standard approaches to semiparametric inference.
In cases with selection on observables we do not need an instrument (or can simply set $Z_i = W_i$),
so for simplicity of notation we replace all instances of $Z_i$ with $W_i$.

\paragraph{Binary treatment with selection on observables.} Most existing work on policy learning,
including \citet{kitagawa2015should}, has focused on the setup where $W_i$ is binary and unconfounded,
i.e., $\cb{Y_i(0), \, Y_i(1)} \indep W_i \cond X_i$. In this case, weighting by the inverse
propensity score lets us recover the average treatment effect, i.e., $g(x, \, w) = (w - e(x)) / (e(x)(1 - e(x)))$
with $e(x) = \PP{W_i = 1 \cond X_i = x}$ identifies the conditional average treatment effect
\smash{$\tau_m(x) = m(x, \, 1) - m(x, \, 0)$} via \eqref{eq:representer}. The estimation
strategy \eqref{eq:LR} yields
\begin{equation}
\htheta = \frac{1}{n} \sum_{i = 1}^n  \p{\hatm(X_i, \, 1) - \hatm(X_i, \, 0) + \frac{W_i - \he(X_i)}{\he(X_i)(1 - \he(X_i)}\p{Y_i - \hatm\p{X_i, \, W_i}}},
\end{equation}
and recovers augmented inverse propensity weighting \citep*{robins1994estimation}.

\paragraph{Continuous treatment with selection on observables.}
In the case where $W_i$ is continuous and unconfounded $\cb{Y_i(w)} \indep W_i \cond X_i$,
we can derive a representer $g(\cdot)$ via integration by parts \citep*{powell1989semiparametric}.
Under regularity conditions, the $\tau$-function
\smash{$\tau_m(x, \, w) = [{d}/{d\nu} \ m(x, \, w + \nu)]_{\nu = 0}$} can be identified
via \eqref{eq:representer} using
\begin{equation}
\label{eq:continuous_riesz}
\begin{split}
&\int \int \frac{d}{dw} \sqb{m(X_i, \, W_i)}_{w = W} dF_{W_i | X_i} \ dF_{X_i}
=  \int \int  g(X_i, \, W_i) m(X_i, \, W_i) dF_{W_i| X_i} \ dF_{X_i}, \\
&g(X_i, \, W_i) = -\frac{d}{dw} \sqb{\log\p{f\p{w \cond X_i}}}_{w = W_i},
\end{split}
\end{equation}
where $f(\cdot \cond x)$ denotes the conditional density of $W_i$ given $X_i = x$.
The resulting doubly robust estimator was to our knowledge first derived via the general approach of \citet{chernozhukov2016locally}, which in turn is closely related to an approach proposed by \citet{ai2007estimation}.

\paragraph{Binary, endogenous treatment with binary treatment and instrument.}
Instead of unconfoundedness, now suppose that $Z_i$ is a valid instrument conditionally on
features $X_i$ in the sense of Assumption 2.1 of \citet{abadie2003semiparametric}. Suppose
moreover that treatment effects are homogenous, meaning that the conditional average treatment
effect matches the conditional local average treatment effect \citep{imbens1994late},\footnote{As
discussed above, our notation has potential outcomes $Y_i(W_i)$ that only depend on treatment $W_i$, and
do not involve the instrument $Z_i$. This is only meaningful when the exclusion restriction holds.}
\begin{equation}
\label{eq:iv_homog}
\tau_m(x) =  m(x, \, 1) - m(x, \, 0) = \frac{\Cov{Y_i, \, Z_i \cond X_i = x}}{\Cov{W_i, \, Z_i \cond X_i = x}}.
\end{equation}
Then we can use a weighting 
function $g(\cdot)$ defined in terms of the compliance score \citep{abadie2003semiparametric,aronow2013beyond},
\begin{equation}
\label{eq:iv_weight}
\begin{split}
&g(X_i, \, Z_i) = \frac{1}{\Delta(X_i)} \frac{Z_i - z(X_i)}{z(X_i)(1 - z(X_i)}, \ \ z(x) = \PP{Z_i = 1 \cond X_i = x}, \\
&\Delta(x) = \PP{W_i = 1 \cond Z_i = 1, \, X_i = x} -  \PP{W_i = 1 \cond Z_i = 0, \, X_i = x},
\end{split}
\end{equation}
to identify this $\tau$-function using \eqref{eq:representer}. We note that our formal results all require that
$g(\cdot)$ be bounded, which implicitly rules out the case of weak instruments (since if $\Delta$ approaches
0, the $g(\cdot)$-weights blow up).

\subsection{Assumptions about the Policy Class}
\label{sec:donsker}

Next, in order to obtain regret bounds that decay as $1/\sqrt{n}$, we need some
control over the complexity of the class $\Pi$ (and again let $\Pi$ potentially change with $n$
for generality). The Vapnik-Chervonenkis (VC) approach \citep{vapnik2000nature} presents us with a
natural way to do so. Recall that the VC-dimension of a class $\Pi$ of binary decision rules
is the largest value of $d \in \NN$ such that there exists a set of $d$ points $x_1, \, ..., \, x_d \in \xx$
that is ``shattered'' by $\Pi$ in the following sense: For each $2^d$ of the binary vectors $v \in \cb{0, \, 1}^d$,
there exists a policy $\pi_v \in \Pi$ such that $\pi_v(X_i) = v_i$ for all $i = 1, \, ..., \, d$.
Throughout our analysis, we control the complexity of $\Pi_n$ by assuming that its VC-dimension
does not grow too fast with the sample size $n$.
As is familiar from the literature on classification, we will
find that the best possible uniform regret bounds scale as $\sqrt{\VC(\Pi_n)/n}$ \citep{vapnik2000nature}.

\begin{assu}
\label{assu:VC}
We assume that there are constants $ 0 < \beta < 1/2$ and $N \geq 1$ such that the Vapnik-Chervonenkis dimension of $\Pi_n$ is
bounded as $\VC(\Pi_n) \leq n^\beta$ for all $n \geq N$.
\end{assu}

In order to illustrate this assumption, we give two examples of policy classes that have a finite VC dimension, and
one that does not. In all three examples below, we assume that the features $X_i$ take values in
$\xx = \RR^p$ for some $p \geq 1$.

\paragraph{Linear Rules}
The VC-dimension of the class of linear decision rules is \citep[][p.~116]{wainwright2019high}
$\VC(\Pi) = p+1$ for
 $\Pi = \cb{\pi_{v,c} :  \pi_{v,c}(x) =  1\p{\cb{v \cdot x \geq c}}, \, v \in \RR^p, \, c \in \RR}$.
Thus, our approach applies to linear decision rules in dimension $p_n \leq n^\beta$ for some $\beta < 1/2$.

\paragraph{Decision Trees}

Trees represent decision rules recursively \citep*{breiman1984classification}. A depth-0 decision
tree $T_0$ is a trivial decision rule, $T_0(x) = a$ for some $a \in \cb{0, \, 1}$ and all $x \in \xx$. For
any $L \geq 1$, a depth-$L$ decision tree $T_L$ is specified via a splitting variable $j \in 1, \, ..., \, p$,
a threshold $t \in \RR$, and two depth-$(L-1)$ decision trees $T_{(L-1),A}$ and $T_{(L-1),B}$, such that
$T_L(x) = T_{(L-1),A}(x)$ if $x_j \leq t$, and $T(x) = T_{(L-1),B}(x)$ else. See Figure \ref{fig:tree} for
an example of a decision tree. The class of depth-$L$ decision trees over $\RR^p$ has VC dimension
bounded on the order of $\VC(\Pi) = \too\p{2^L \log(p)}$.\footnote{This bound follows Lemma 4 of
\citet*{zhou2018offline}, paired with the alternative characterization of the VC dimension given in
Section A of the supplemental material. The notation $f(n) = \too(g(n))$ means that there is a  function $h(\cdot)$
that scales poly-logarithmically in its argument for which $f(n) \leq h(g(n))g(n)$.}
Thus, our results apply to trees whose depth may grow
as $L_n = \lfloor \kappa \log_2(n) \rfloor$ for some $\kappa < 1/2$.

\paragraph{Monotone Rules}

We have $x \in [0, \, 1]^2$ and units get treated if $x_2$
exceeds some increasing function of $x_1$, i.e.,
$\Pi = \cb{\pi_f : \pi_f(x) = 1\p{\cb{x_2 \geq f(x_1)}}, \, \text{$f$ is monotone increasing}}$.
This class has infinite VC dimension, because any set of points $\cb{x_i}_{i = 1}^d$ with
$x_i = (\alpha_i, \, \alpha_i^2)$ and $0 < \alpha_1 < \ldots < \alpha_d < 1$ can be shattered
using $\Pi$. Thus, our results do not apply to monotone rules over $[0, \, 1]^2$.\footnote{The 
difficulty here is not a mere technicality: Monotone decision rules can match
arbitrary decision rules along the curve $(\alpha, \, \alpha^2)$ for $\alpha \in [0, \, 1]$,
and so it is impossible to establish any non-trivial learning rates over monotone decision rules
without making further assumptions on the distribution of the features $X_i$. In particular, we need
assumptions that guarantee that all observations cannot concentrate around the curve $(\alpha, \, \alpha^2)$.
In this paper, we do not consider results that require specific distributional assumptions over the
features $X_i$. We note however the recent work by \citet{mbakop2016model}, who establish polynomial
rates of convergence for learning monotone rules under an assumption that the $X_i$ have a bounded
density under Lebesgue measure on $[0, \, 1]^2$.}

\subsection{Bounding Asymptotic Regret}
\label{sec:main_result}

We are now ready to state our main result on the asymptotic regret of policy learning
using doubly robust scores. Following \citet{chernozhukov2016double,chernozhukov2016locally}
we assume that we run our method with scores obtained via cross-fitting, which is a type of data
splitting that can be used to verify asymptotic normality given only high-level conditions on the
predictive accuracy of the methods used to estimate nuisance components.
In particular, cross-fitting allows for the use of black-box machine learning tools provided we can verify
that they are accurate in mean-squared error as in Assumption \ref{assu:LR}.

We proceed as follows:
First divide the data into $K$ evenly-sized folds and, for each fold $k = 1, \, ..., \, K$,
run an estimator of our choice on the other $K - 1$ data folds to estimate the
functions \smash{$m_n(x, \, w)$} and \smash{$g_n(x, \, z)$}; denote the resulting estimates 
\smash{$\hatm_n^{(-k)}(x, \, w)$} and \smash{$\hg_n^{(-k)}(x, \, z)$}.  Throughout, we will only
assume that these nuisance estimates are accurate in the sense of Assumption \ref{assu:LR}.
Then, given these pre-computed values, we choose \smash{$\hpi_n$} by maximizing a doubly
robust estimate of $A(\pi) = 2V(\pi) - \EE{\tau(X_i)}$,
\begin{equation}
\label{eq:estimator}
\begin{split}
&\hpi_n = \argmax\cb{\hA_n(\pi) : \pi \in \Pi_n}, \ \ \hA_n(\pi) = \frac{1}{n} \sum_{i = 1}^n \p{2\pi(X_i) - 1} \hGamma_i, \\
&\hGamma_i = \tau_{\hatm_n^{(-k(i))}}(X_i, \, W_i) + \hg_n^{(-k(i))}\p{X_i, \, Z_i} \p{Y_i - \hatm_n^{(-k(i))}\p{X_i, \, W_i}},
\end{split}
\end{equation}
where $k(i) \in \cb{1, \, ..., \, K}$ denotes the fold containing the $i$-th observation.
The $K$-fold algorithmic structure used in
\eqref{eq:estimator} was proposed in an early paper by \citet{schick1986asymptotically} as a general purpose
tool for efficient estimation in semiparametric models, and has also been used by other authors including
\citet{robins2017minimax} and \citet{zheng2011cross}.

Finally, we assume that the weighting function $g_n(x, \, z)$ is bounded uniformly as below.
In the case of a binary exogenous treatment, this is equivalent to the ``overlap'' assumption
in the causal inference literature \citep{imbens2015causal}, whereby $\eta \leq \PP{W_i = 1 \cond X_i = x} \leq 1 - \eta$
for all values of $x$. In our setting, the condition below
acts as a generalization of the overlap assumption \citep{hirshberg2018balancing}.

\begin{assu}
\label{assu:overlap}
There is an $\eta > 0$ such that
$\abs{g_n(x, \, z)} \leq \eta^{-1} \text{ for all } x, \, z, \, n$.
\end{assu}

We also define the following quantities, where
$S_n$ bounds the second moment of the scores, and $S_n^*$ is the asymptotic
variance for estimating the policy improvement $A(\pi)$ of the best policy
in $\Pi_n$ via \eqref{eq:LR}:\footnote{By expanding the square, we see that policies with higher values have
lower variance of their scores, and so $S_n^*$ corresponds to the asymptotic variance for
evaluating an optimal policy. Moreover, in the case where arguments from \citet{newey1994asymptotic} imply
that the doubly robust estimator \eqref{eq:LR} is efficient, then $S_n^*$ is the semiparametric
efficient variance for evaluating an optimal policy.}
\begin{align}
\label{eq:Sn}
&S_n = \EE{\p{\tau_{m_n}(X_i, \, W_i) - g_n(X_i, \, Z_i)\p{Y_i - m_n(X_i, \, W_i)}}^2}, \\
\notag
&S_n^* = \inf\cb{\Var{(2\pi(X_i) - 1) \p{\tau_{m_n}(X_i, \, W_i) - g_n(X_i, \, Z_i)\p{Y_i - m_n(X_i, \, W_i)}}} : \pi \in \Pi_n}.
\end{align}
We note that, unless we have an exceptionally large signal-to-noise ratio, we will have
$S_n^* \geq S_n / 4$ and so the rounded log-term in \eqref{eq:theo_main} below is just 0.
A proof of Theorem \ref{theo:main} is given in the following section.

\begin{theo}
\label{theo:main}
Given Assumptions \ref{assu:identification},  \ref{assu:LR} and \ref{assu:overlap},
define \smash{$\hpi_n$} as in \eqref{eq:estimator}.\footnote{We assume
that the rates of convergence specified in Assumption \ref{assu:LR} apply to the nuisance
components estimated for each fold $k = 1, \, ..., \, K$ in \eqref{eq:estimator}.}
Suppose moreover that the irreducible noise $\varepsilon_i = Y_i -  m(X_i, \, W_i)$
is both uniformly sub-Gaussian conditionally on $X_i$ and $W_i$ and has second moments
uniformly bounded from below, $\Var{\varepsilon_i \cond X_i=x, \, W_i=w} \geq s^2$, and that
the treatment effect function \smash{$\tau_{m_n}(x, \, w)$} is uniformly bounded in $x$, $w$ and $n$.
Finally, suppose that $\Pi_n$ satisfies Assumption \ref{assu:VC} with parameter
$\beta \leq \min\cb{\zeta_m, \, \zeta_g}$, where the $\zeta$ are as defined in Assumption \ref{assu:LR}.
Then, for any sequence $\psi_n \geq 0$ with $\limn \psi_n \sqrt{n} = 0$,
\begin{equation}
\label{eq:theo_main}
\begin{split}
&\limsup_{n \rightarrow \infty} \ \EE{\sup\cb{R_n\p{\pi} : \hA_n(\pi) \geq \max\cb{\hA_n(\pi) : \pi \in \Pi_n} - \psi_n, \ \ \pi \in \Pi_n}} \\
&\ \ \ \ \ \ \ \ \ \ \ \Bigg/ \sqrt{\VC(\Pi_n) S_n^* \p{1 + \left\lfloor \log_4\p{\frac{S_n}{S_n^*}} \right \rfloor \Big/\, 9} \Big/\, n} \leq 60,
\end{split}
\end{equation}
where $R_n(\cdot)$ denotes regret for the $n$-th data-generating distribution.
\end{theo}

In the simplest case where the maximizer of \smash{$\hA_n(\pi)$} over $\pi \in \Pi_n$
is unique and $\psi_n = 0$ (i.e., we solve the maximization problem exactly), the statement
in \eqref{eq:theo_main} simplifies to a bound on \smash{$\EE{R_n\p{\hpi_n}}$},
where \smash{$\hpi_n$} is as defined in \eqref{eq:estimator}. However, in practice, \smash{$\hA_n(\pi)$}
may have many maximizers. Moreover, the optimization problem \eqref{eq:estimator} is not convex and
so---given a reasonable computational budget---we may only be able to solve it to within some
tolerance $\psi_n > 0$. The more comprehensive form of our result given above highlights the fact that,
in this case, our regret bound in fact applies uniformly over all approximate solutions to \eqref{eq:estimator}.

\section{Upper Bounds}
\label{sec:formalresults}

In this section, we present a series of results that culminate in a proof of Theorem \ref{theo:main},
given in Section \ref{sec:theo1}. All other proofs are deferred to Section C of the supplemental material.
Recall that we study policy learning for a class of problems where regret can be
written as in \eqref{eq:regret} using a function $V_n(\pi) = \EE[n]{\pi(X_i)\tau_n(X_i)}$,
and we obtain \smash{$\hpi_n$} by maximizing a cross-fitted doubly robust estimate of
\smash{$A_n(\pi) = 2V_n(\pi) - \EE[n]{\tau_n(X_i)}$} defined in \eqref{eq:estimator} over the class $\Pi_n$.
If we could use $\hA_n(\pi) = A_n(\pi)$, then \eqref{eq:estimator} would directly yield the regret-minimizing policy
in the class $\Pi_n$; but of course we never know $A_n(\pi)$ in applications. Thus, the main focus of our formal
results is to study stochastic fluctuations of the empirical process \smash{$\hA_n(\pi) - A_n(\pi)$} for $\pi \in \Pi_n$,
and examine how they affect the quality of policies learned via \eqref{eq:estimator}.

\subsection{Rademacher Complexities and Oracle Regret Bounds}
\label{sec:influence}

We start our analysis by characterizing concentration of
an ideal version of the objective in \eqref{eq:estimator} based on the true
influence scores $\Gamma_i$, rather than doubly robust estimates thereof:
\begin{equation}
\label{eq:influence}
\tA_n(\pi) = \frac{1}{n} \sum_{i = 1}^n \p{2\pi(X_i) - 1} \Gamma_i, \ \ 
\Gamma_i = \tau_{m_n}(X_i, \, W_i) + g_n\p{X_i, \, Z_i} \p{Y_i - m_n\p{X_i, \, W_i}}.
\end{equation}
The advantage of studying concentration of the empirical process
\smash{$\tA_n(\pi) - A_n(\pi)$} over the set $\pi \in \Pi_n$ is that it allows us, for the time being,
 to abstract away from the estimation tools used to obtain \smash{$\hA_n(\pi)$},
 and instead to focus on the complexity of empirical maximization over the class $\Pi_n$.

A convenient way to bound the supremum of this
empirical process over any class $\Pi$ is by controlling its Rademacher complexity
$\rr_n(\Pi)$, defined as\footnote{Note that, conditionally on $\cb{X_i, \, \Gamma_i}_{i = 1}^n$
and the Rademacher variables $\xi_i$, the sum $\sum_{i = 1}^n \xi_i \Gamma_i \p{2\pi(X_i) - 1}$
can only take $2^n$ distinct values. Thus, the definition of $\rr_n(\Pi)$ does not entail any
measure theoretic problems.}
\begin{equation}
\label{eq:rademacher}
\rr_n(\Pi) =  \EE{\sup_{\pi \in \Pi} \cb{\frac{1}{n} \sum_{i = 1}^n \xi_i \Gamma_i \p{2\pi(X_i) - 1}} \cond \cb{X_i, \, \Gamma_i}_{i = 1}^n}
\end{equation}
where the $\xi_i$ are independent Rademacher (i.e., sign) random variables $\xi_i = \pm 1$
with probability $1/2$ each \citep{bartlett2002rademacher}.
For intuition as to why Rademacher complexity is a natural complexity measure,
note that $\rr_n(\Pi)$ characterizes the maximum (weighted) in-sample
classification accuracy on randomly generated labels $\xi_i$ over classifiers $\pi \in \Pi$;
thus, $\rr_n(\Pi)$ measures how much we can overfit to random coin flips using $\Pi$.

Following this proof strategy, we bound the Rademacher complexity of ``slices''
of our policy class $\Pi_n$, defined as
\begin{equation}
\label{eq:pi_lambda}
\Pi^\lambda_n = \cb{\pi \in \Pi_n : R_n\p{\pi} \leq \lambda}.
\end{equation} 
The reason we focus on slices of $\Pi_n$ is that, when we use doubly robust scores,
low-regret policies can generally be evaluated more accurately than high-regret policies,
and using this fact allows for sharper bounds. Specifically, we can check that
\smash{$n\text{Var}[\tA_n(\pi)] = S_n - A_n^2(\pi)$}, and so
\begin{equation}
\label{eq:SL}
n\sup\cb{\Var{\tA_n(\pi)} : \pi \in \Pi_n^\lambda} := S_n^\lambda \leq S_n^* + 4\lambda \sup\cb{A_n(\pi) : \pi \in \Pi_n},
\end{equation}
where $S_n$ and $S_n^*$ are defined in \eqref{eq:Sn}.
This type of slicing technique is common in the literature, and has been used in different
contexts by, e.g., \citet*{bartlett2005local} and \citet{gine2006concentration}.

The following result provides such a bound in terms of the second moments of the doubly robust score,
specifically $S_n^\lambda$ and $S_n$. This bound is substantially stronger than corresponding
bounds used in existing results on policy learning. \citet{kitagawa2015should} build
their result on bounds that depend
on $\max\cb{\Gamma_i}/\sqrt{n}$, which can only be used with scores that are uniformly 
bounded in order to get optimal rates.
Meanwhile, bounds that scale as $\sqrt{S_n^\lambda\log(n)/n}$ are developed
by \citet*{cortes2010learning}, \citet{maurer2009empirical} and \citet{swaminathan2015batch}; however,
the additional $\log(n)$ factor makes these bounds inappropriate for asymptotic analysis.

\begin{lemm}
\label{lemm:rademacher}
Suppose that the class $\Pi_n$ satisfies Assumption \ref{assu:VC}, and that 
the scores \smash{$\Gamma_i$} in \eqref{eq:influence} are drawn from a sequence
of uniformly sub-Gaussian distributions with variance bounded from below,
\begin{equation}
\label{eq:subG}
\PP[n]{\abs{\Gamma_i} > t} \leq C_\nu  \, e^{-\nu t^2} \text{ for all } t > 0, \ \ \Var[n]{\Gamma_i \cond X_i = x} \geq s^2,
\end{equation}
for some constants $C_\nu, \, \nu, \, s > 0$ and all $n= 1, \, 2, \, ...$
Then, for any $\lambda$,
\begin{equation}
\label{eq:rad_bound}
\limsup_{n \rightarrow \infty} \ \EE{\rr_n\p{\Pi_n^\lambda}}  \, \bigg/\, \sqrt{\p{S_n^\lambda + 4\lambda^2} \p{1 + \left\lfloor\log_4\p{ \frac{S_n}{S_n^\lambda}} \right\rfloor \Big/\,9} \frac{\VC(\Pi_n)}{n}}  \, \leq \, 20.
\end{equation} 
\end{lemm}

Then, following the well known approach of \citet{bartlett2002rademacher}, we use our bound on
Rademacher complexity to obtain a uniform concentration bound for \smash{$\tA_n(\pi)$}.
We use a refinement of the argument of \citet{bartlett2002rademacher} based on Talagrand's inequality to
get a bound that depends on second moments of $\Gamma_i$ rather than $\sup \abs{\Gamma_i}$.

\begin{coro}
\label{coro:additive}
Under the conditions of Lemma \ref{lemm:rademacher}, the expected maximum error of
\smash{$\tA_n(\pi)$} is bounded as
\begin{equation}
\label{eq:influence_concentration}
\begin{split}
&\limsup_{n \rightarrow \infty} \ \EE{\sup \cb{ \abs{\tA_n\p{\pi} - A_n\p{\pi}} : \pi \in \Pi_n^\lambda}}  \\
& \ \ \ \ \ \ \ \ \ \ \ \ \ \ \ \ \ \ \ \ \ \ \ \ 
\Bigg/ \sqrt{\p{S_n^\lambda + 4\lambda^2}  \p{1 + \left\lfloor\log_4\p{ \frac{S_n}{S_n^\lambda}} \right\rfloor \Big/\, 9} \frac{\VC(\Pi_n)}{n}}  \, \leq \, 40.
\end{split}
\end{equation} 
Furthermore, this error is concentrated around its expectation:
There is a sequence $c_n \rightarrow 0$ such that, for any $\delta > 0$,
\begin{equation}
\label{eq:influence_concentration_sharp}
\begin{split}
&\sup \cb{ \abs{\tA_n\p{\pi} - A_n\p{\pi}} : \pi \in \Pi_n^\lambda} \\
& \ \ \ \ \ \ \ \ \ 
\leq \p{1 + c_n}\p{\EE{\sup \cb{ \abs{\tA_n\p{\pi} - A_n\p{\pi}} : \pi \in \Pi_n^\lambda}} + \sqrt{\frac{2 S_n^\lambda \log(\delta^{-1})}{n}}}
\end{split}
\end{equation} 
with probability at least $1 - \delta$.
\end{coro}

In our final argument, we will apply Corollary \ref{coro:additive} for different $\lambda$-slices,
and verify that we can in fact focus on those slices where $\lambda$ is nearly 0. Before that,
however, we also need to control the discrepancy between the feasible objective \smash{$\hA_n(\pi)$}
and the oracle surrogate \smash{$\tA_n(\pi)$} studied here.

\subsection{Uniform Coupling with the Doubly Robust Estimator}
\label{sec:coupling}

In the previous section, we established risk bounds that would hold if we could
optimize the infeasible value function \smash{$\tA_n(\pi)$}; we next need to extend
these bounds to cover the situation where we optimize a feasible value function. As
discussed above, we focus on the doubly robust estimator \eqref{eq:estimator}, obtained
using cross-fitting as in \citet{chernozhukov2016double,chernozhukov2016locally}.
As preliminaries, we note that the results of \citet{chernozhukov2016locally}
immediately imply that, given Assumption \ref{assu:LR}, \smash{$\hA_n(1)$} is an asymptotically
normal estimate of $A_n(1)$, where we use ``1'' as shorthand for the ``always treat'' policy.
 Furthermore, it is easy to check that given any fixed policy $\pi$,
\begin{align}
\label{eq:coupling}
\sqrt{n} \p{\hA_n\p{\pi} - \tA_n\p{\pi}} \rightarrow_p 0,
\end{align}
meaning that the discrepancy between the two value estimates decays faster than
the variance of either.

However, in our setting, the analyst gets to optimize over all policies $\pi \in \Pi_n$,
and so coupling results established for a single pre-determined policy $\pi$ are not
strong enough.
The following lemma extends the work of \citet{chernozhukov2016locally} to the case
where we seek to establish a coupling of the form \eqref{eq:coupling} that holds simultaneously
for all $\pi \in \Pi_n$.

\begin{lemm}
\label{lemm:coupling}
Under the conditions of Lemma \ref{lemm:rademacher}, suppose that
Assumptions \ref{assu:identification} and \ref{assu:overlap} hold, and that we obtain \smash{$\hA_n(\pi)$} using
cross-fitted estimates of nuisance components satisfying Assumption \ref{assu:LR}. Then
\begin{equation}
\label{eq:unif_coupling}
\begin{split}
\frac{\sqrt{n} \  \EE{\sup\cb{\abs{\hA_n(\pi) - \tA_n\p{\pi}} : \pi \in \Pi_n}}}{ a\p{\p{1 - K^{-1}} n}} = \oo\p{1 + \sqrt{\frac{\VC(\Pi_n)}{n^{\min\cb{\zeta_m, \, \zeta_g}}}}},
\end{split}
\end{equation}
where the $\oo(\cdot)$ term hides a dependence on the overlap parameter $\eta$ from Assumption \ref{assu:overlap}
and the sub-Gaussianity parameter $\nu$ specified in Lemma \ref{lemm:rademacher}.
\end{lemm}

The above result is perhaps surprisingly strong: Provided that the dimension
$\VC(\Pi_n)$ of $\Pi_n$ does not grow too fast with $n$,
the bound \eqref{eq:unif_coupling} is the same coupling bound as we might expect to obtain for a single
policy $\pi$, and the dimension of the class $\Pi_n$ does
not affect the leading-order constants in the bound. In other words, in terms of the coupling of
\smash{$\tA_n(\pi)$} and \smash{$\hA_n(\pi)$}, we do not lose anything by scanning over a
continuum of policies $\pi \in \Pi_n$ rather than just considering a single policy $\pi$.

The doubly robust form used here
is not the only way to construct efficient estimators for the value of a single policy $\pi$---for example,
\citet*{hirano2003efficient} show that inverse-propensity weighting with non-parametrically estimated
propensity scores may also be efficient---but it plays a key role in the proof of Lemma \ref{lemm:coupling}.
In particular, under Assumption \ref{assu:LR}, the natural bound for the bias term due to misspecification
of the nuisance components in fact holds simultaneously for all $\pi \in \Pi$, and this helps us pay a
smaller-than-expected price for seeking a uniform result as in \eqref{eq:unif_coupling}.
It is far from obvious that other efficient methods
for evaluating a single policy $\pi$, such as that of \citet{hirano2003efficient}, would lead to equally
strong uniform couplings over the whole class $\Pi_n$.

\subsection{Proof of Theorem 1}
\label{sec:theo1}

Given that Assumption \ref{assu:identification}, \ref{assu:LR}, \ref{assu:VC} and \ref{assu:overlap} hold with parameters
$\beta < \min\cb{\zeta_m, \, \zeta_g}$, a combination of results from  Corollary \ref{coro:additive}
and Lemma \ref{lemm:coupling} implies that \smash{$\hA_n(\cdot)$} concentrates around $A_n(\cdot)$
over $\Pi_n^\lambda$.
To conclude, it now remains to apply these bounds at two different values of $\lambda$.
First we choose $\lambda^* > 0$ such as to satisfy
$4(\lambda^*)^2 + 4\lambda^*\sup\cb{A(\pi) : \pi \in \Pi_n} \leq S_n^*$, so that the following
holds via \eqref{eq:SL}:
$$ S_n^{\lambda^*} + 4(\lambda^*)^2 \leq S_n^* + 4(\lambda^*)^2 + 4\lambda^*\sup\cb{A(\pi) : \pi \in \Pi_n} \leq 2S_n^*. $$
Then, by Corollary \ref{coro:additive} and Lemma \ref{lemm:coupling}, we find that the limsup of the
following expression is bounded by 1 as $n$ goes to infinity:
\begin{equation*}
 \EE{\sup \cb{ \abs{\hA_n\p{\pi} - A_n\p{\pi}} : \pi \in \Pi_n^{\lambda^*}}} \Bigg/ \p{60 \sqrt{S_n^* \p{1 + \left\lfloor\log_4\p{ \frac{S_n}{S_n^*}} \right\rfloor \Big/\, 9} \frac{\VC(\Pi_n)}{n}}}.
\end{equation*} 
Now, recall that if any two functions $h(\cdot)$ and \smash{$\hat{h}(\cdot)$} are uniformly coupled as
\smash{$|h(u) - \hat{h}(u)| \leq b$} for all $u \in U$ and \smash{$\hat{h}(\hat{u}) \geq \sup\{\hat{h}(u) : u \in U\} - \psi$}, then
$$h(\hat{u}) \geq \hat{h}(\hat{u}) - b \geq \hat{h}(u) - b - \psi \geq h(u) - 2b - \psi$$
for any $u \in U$. Thus, the above implies that (recall that $A_n(\pi)$ scales with $2R_n(\pi)$)
\begin{equation}
\label{eq:abcd}
\begin{split}
&\limsup_{n \rightarrow \infty} \ \EE{\sup\cb{R_n\p{\pi} : \hA_n(\pi) \geq \max\cb{\hA_n(\pi) : \pi \in \Pi_n^{\lambda^*}} - \psi_n, \ \ \pi \in \Pi_n^{\lambda^*}}} \\
&\ \ \ \ \ \ \ \ \ \ \ \Bigg/ \p{\frac{\psi_n}{2} +  60 \sqrt{S_n^* \p{1 + \left\lfloor\log_4\p{ \frac{S_n}{S_n^*}} \right\rfloor \Big/\, 9} \frac{\VC(\Pi_n)}{n}}} \leq 1,
\end{split}
\end{equation} 
and we note that $\psi_n$ decays fast enough by assumption that it can be omitted from \eqref{eq:abcd} without altering the result.
In other words, if we knew that our learned policy approximately maximizes \smash{$\hA_n(\pi)$} \emph{and}
has regret less than $\lambda^*$, then we could guarantee that its regret decays at the desired rate.

To prove our result, it remains to show that all approximate maximizers of $\hA_n(\cdot)$ have regret bounded by $\lambda^*$
enough for \eqref{eq:abcd} to capture the leading-order behavior of regret. To do so, we apply
a similar argument as above, but at a different value of $\lambda$. Consider
\smash{$\lambda_+ = 3\limsup_{n \rightarrow \infty} \sup\cb{R_n(\pi) : \pi \in \Pi_n}$},
and by \eqref{eq:influence_concentration_sharp} we see that
\begin{equation}
\label{eq:lambdaplus}
\limn \sqrt{n} \, \PP{\sup\cb{\abs{\tA_n(\pi) - A_n(\pi)} : \pi \in \Pi_n^{\lambda_+}} \geq \frac{\lambda^*}{5}} = 0.
\end{equation}
Furthermore, note that \smash{$\Pi_n^{\lambda_+} = \Pi_n$} for large enough $n$, and so \eqref{eq:lambdaplus}
in fact also holds with  \smash{$\Pi_n^{\lambda_+}$} replaced by $\Pi_n$.
Meanwhile, from \eqref{eq:unif_coupling} paired with Markov's inequality we know that
\begin{equation}
\PP{\sup\cb{\abs{\hA_n(\pi) - \tA_n(\pi)} : \pi \in \Pi_n} \geq\frac{\lambda^*}{5}} = \oo\p{\frac{a\p{\p{1 - K^{-1}}n}}{\sqrt{n}}}.
\end{equation}
By combining these two bounds, we see that
\begin{equation}
\label{eq:better_than_lP}
\begin{split}
&\limn \sqrt{n} \, \mathbb{P}\bigg[ \cb{\pi \in \Pi_n : \hA_n(\pi) \geq \max\cb{\hA_n(\pi) : \pi \in \Pi_n} - \psi_n} \\
&\ \ \ \ \ \ \ \ \ \ \ \ \ \ \ \ \ \ \ \ \ \  \bigcap \cb{\pi \in \Pi_n : R_n(\pi) \geq \lambda^*} \neq \emptyset\bigg] =  0,
\end{split}
\end{equation}
and moreover,  because $\tau_{m_n}(x, \, w)$ is uniformly bounded, we find that
the contribution of events where \eqref{eq:better_than_lP} fails to hold to \eqref{eq:theo_main}
is vanishingly small as $n$ gets large.

\section{Lower Bounds}
\label{sec:lb}

To complement the upper bounds given in Theorem \ref{theo:main},
we also present lower bounds on the minimax risk for policy learning.
Our goal is to show that our bounds are the best possible regret bounds
that flexibly account for the distribution of the observed data and
depend on the policy class $\Pi$ through the Vapnik-Chervonenkis dimension $\VC(\Pi)$.
For simplicity, we here only consider the case where $W_i$ is binary and unconfounded; lower
bounds for other cases considered in this paper can be derived via analogous arguments.

To establish our result, we consider lower bounds over sequences of problems defined as follows.
Let $\xx_s := [0, \, 1]^s$ denote the $s$-dimensional unit cube for some positive integer $s$, and let
$f(x)$ and $e(x)$ be $\lceil s/2 + 1 \rceil$ times continuously differentiable
functions over $\xx_s$. Moreover, let $\sigma^2(x)$ and $\tau(x)$ be functions on $\xx_s$
such that $\sigma^2(x)$ is bounded away from 0 and $\infty$, and
$\abs{\tau(x)}$ is bounded away from $\infty$. Then, we define an
asymptotically ambiguous problem sequence as one where $\cb{X_i, \, Y_i, \, W_i}$
are independently and identically distributed drawn as
\begin{equation}
\label{eq:aa}
\begin{split}
&X_i \sim \pp, \ \ W_i \cond X_i \sim \text{Bernoulli}(e(X_i)), \\
&Y_i \cond X_i, \, W_i \sim \nn\p{f(X_i) + \p{W_i - e(X_i)} \frac{\tau(X_i)}{\sqrt{n}}, \, \sigma^2(X_i)}.
\end{split}
\end{equation}
Because of the number of derivatives assumed on $f(x)$ and $e(x)$, it is well known
that simple series estimators satisfy Assumption \ref{assu:LR}.\footnote{See \citet{nickl2007bracketing}
for an argument that holds for arbitrary distributions $\pp$ supported on $[0, \, 1]^s$. We also note that,
for a complete argument, one needs to address the fact that we have not
assumed the treatment effect function $\tau(x)$ to be differentiable.
To address this issue, note that in our data-generating process \eqref{eq:aa} we have $\EE{Y_i|X_i= x} = f(x)$ regardless
of $n$. Thus, because both $e(x)$ and $f(x)$ are sufficiently differentiable, we can use
standard results about series estimation to obtain $o_P(n^{-1/4})$-consistent
estimators \smash{$\he(x)$} and \smash{$\hat{f}(x)$} for these quantities. Next, for the purpose of
our policy learner, we simply set \smash{$\hatm(x, 0) = \hatm(x, 1) = \hat{f}(x)$}; and because
$\EE{\tau^2(X_i)/\sqrt{n}} = \oo(1/n)$, these regression adjustments in fact
satisfy Assumption \ref{assu:LR}.}
Thus, because the magnitude of the treatment effects shrinks in \eqref{eq:aa},
\smash{$S_n^*$} and $S_n$  both converge to $S_{\pp}$ as defined below, and so Theorem \ref{theo:main}
immediately implies that, under unconfoundedness,
\begin{equation}
\label{eq:aabound}
\limsup_{n \rightarrow \infty} \ R_n\p{\hpi_n} \, \bigg/ \sqrt{\frac{S_{\pp} \VC\p{\Pi}}{n}} \leq 60, \ \ S_{\pp} = \EE[\pp]{\frac{\sigma^2(X_i)}{e(X_i)\p{1 - e(X_i)}}}
\end{equation}
for any policy class $\Pi$ with finite VC dimension. 
The following result shows that \eqref{eq:aabound} is sharp up to a
universal constant (whose value is less than 200).\footnote{The strategy of proving lower bounds relative to an adversarial
feature distribution $\pp$ is standard in the machine learning literature; see, e.g.,
\citet{devroye1995lower}. If we fix the distribution $\pp$ a-priori, then regret
bounds for empirical risk minimization over $\Pi$ based on structural
summaries of $\Pi$ (such as the VC dimension) may be loose \citep{bartlett2006empirical};
however, it is not clear how to exploit this fact other than by conducting ad-hoc analyses
for specific choices of $\Pi$.}

\begin{theo}
\label{theo:lb}
Let $f(x)$, $e(x)$, and $\sigma(x)$ be functions over $\xx_s$ satisfying the conditions
discussed above, and let $\Pi$ be a class of functions  over $\xx_s$ with finite VC dimension.
Then, there exists a distribution $\pp$ supported on $[0, \, 1]^s$ (and a constant $C$)
such that the minimax risk for policy learning over the data generating
distribution \eqref{eq:aa} (with unknown $|\tau(x)| \leq C$) and the policy class $\Pi$
is bounded from below as follows, where $\hpi_n$ can be any measurable function
of the training sample:
\begin{equation}
\label{eq:lb}
 \liminf_{n \rightarrow \infty} \cb{\sqrt{n} \, \inf_{\hpi_n}  \cb{ \sup_{\abs{\tau(x)} \leq C} \cb{\EE{R_n\p{\hpi_n}}}}} \geq 0.33 \sqrt{S_{\pp} \VC\p{\Pi}}.
\end{equation}
\end{theo}

Here, the fact that we focus on problems where
the magnitude of the treatment effect scales as $1/\sqrt{n}$ is important, and closely
mirrors the type of asymptotics used by \citet{hirano2009asymptotics}. If treatment
effects decay faster than $1/\sqrt{n}$, then learning better-than-random policies is effectively
impossible---but this does not matter, because of course all decision rules have regret
decaying as $o(1/\sqrt{n})$ and so Theorem \ref{theo:main} is loose. Conversely, if treatment
effects dominate the  $1/\sqrt{n}$ scale, then in large samples it is all but obvious who
should be treated and who should not, and it is possible to get regret bounds that decay
at superefficient rates \citep{luedtke2017faster}, again making Theorem \ref{theo:main} loose.
But if the treatment effects obey the $\Theta(1/\sqrt{n})$ scaling of \citet{hirano2009asymptotics},
then the problem of learning good policies is neither trivial nor impossible, and the value of using doubly robust policy
evaluation for policy learning becomes apparent.

Finally, we note that the bounds of \citet{kitagawa2015should} for inverse-propensity weighting
are not asymptotically sharp in the above sense. Even when propensity scores are known,
\citet{kitagawa2015should} assume that $\abs{Y_i} \leq M$ and $\eta \leq e(X_i) \leq 1 - \eta$, and then
prove regret bounds that scale as \smash{$M/\eta \sqrt{\VC(\Pi)/n}$} instead of \smash{$\sqrt{S_{\pp}\VC(\Pi)/n}$} in \eqref{eq:aabound}.
Now, the bound of \citet{kitagawa2015should} is of course sometimes sharp, e.g., it
is optimal if all we know is that $\abs{Y_i} \leq M$ and $\eta \leq e(X_i) \leq 1 - \eta$, but it
is not adaptively sharp for asymptotically ambiguous sequences of problems as in \eqref{eq:aa}.
In particular, the ratio of the upper bound of \citet{kitagawa2015should} and the lower bound \eqref{eq:lb}
scales as $M/(\eta \sqrt{S_{\pp}})$, and there exist sequences of type \eqref{eq:aa} where this ratio
may be arbitrarily large.\footnote{Using the
techniques developed in this paper, we can sharpen the bounds of \citet{kitagawa2015should}
and asymptotically replace $M/\eta$ by \smash{$\EE{Y_i^2/(e(X_i)(1 - e(X_i)))}^{1/2}$}. However, even this
improved bound may exceed \eqref{eq:lb} by an arbitrarily large factor.}

\section{Implementation and Experiments}
\label{sec:experiments}

We now illustrate the value of doubly robust scoring techniques for policy
learning using both an example from program evaluation and simulation studies.
In Section \ref{sec:gain} we revisit a randomized evaluation of California's GAIN program, while
Section \ref{sec:simu} presents a simulation study with endogenous treatment assignment.
We present additional simulation results on nudge interventions to a continuous treatment variable
in Section B of the supplemental material.

Recall that our approach to policy involves a 3-step algorithm. We start with
a set of $n$ independent and identically distributed training examples $(X_i, \, Y_i, \, W_i, \, Z_i)$
and a class $\Pi$ of acceptable policies. Then, we
\begin{enumerate}
\item Estimate the nuisance components $m(x, \, w)$ and $g(x, \, z)$ defined in Section \ref{sec:identification},
\item Form doubly robust scores\footnote{Recall that $\tau_m(\cdot)$ does not depend on $w$ in the case of binary
treatments, and we omit the redundant argument in this case.}
$\hGamma_i = \tau_{\hatm}(X_i, \, W_i) + \hg(X_i, \, Z_i) (Y_i - \hatm(X_i, \, W_i))$,
with cross-fitting as discussed in Section \ref{sec:main_result}, and
\item Select $\hpi \in \argmax\cb{\sum_{i = 1}^n (2\pi(X_i)  - 1) \hGamma_i : \pi \in \Pi}$.
\end{enumerate}
The main points of freedom left to the analysts involve the choice of estimator for $m(\cdot)$ and
$g(\cdot)$ in Step 1, and the implementation of the optimization problem in Step 3.
We emphasize that the choice of estimator for $m(x, \, w)$ and $g(x, \, z)$ in Step 1 and
the choice of policy class $\Pi$ along with the optimizer used in Step 3 can be made fully independently.

For Theorem \ref{theo:main} to apply, the main requirement on the method used to estimate
$m(x, \, w)$ and $g(x, \, z)$ in Step 1 is that its error decays fast enough in mean-squared error,
as detailed in Assumption \ref{assu:LR}.
Here, one option is to use non-parametric estimators for which we can precisely spell out when they
satisfy Assumption \ref{assu:LR}, such as sieve-based methods \citep*{chen2007large}
or kernel regression \citep*{caponnetto2007optimal};
another is to use more heuristic methods from the statistical learning literature, such as
boosting, random forests, or neural networks, in the hope that they will empirically be more accurate
in finite samples than sieve or kernel-based methods.\footnote{In a recent
advance, \citet*{farrell2018deep} established conditions under which deep neural networks can be shown to
provably satisfy the conditions required by Assumption \ref{assu:LR}. Thus, depending on the statistical setting and the
chosen architecture, deep neural networks could either be seen as a formally validated alternative to sieve-type methods or as
heuristic method.} One possible compromise is to run both
classical methods known to satisfy Assumption \ref{assu:LR} asymptotically and heuristic
statistical learning tools, and then synthesize the output of all models via cross-validation.
As argued in \citet*{van2007super}, this approach essentially matches the finite-sample accuracy
of the best method under consideration while preserving the asymptotic guarantees of the classical ones.

Meanwhile, the optimization problem in Step 3 is not a convex optimization problem, and so solving
it can be computationally challenging. Several authors, including \citet{beygelzimer2009offset},
\citet{kitagawa2015should}, \citet*{zhang2012estimating} and \citet*{zhao2012estimating},
have noted that this optimization problem is numerically equivalent to a weighted classification problem,
\begin{equation}
\label{eq:classif}
\begin{split}
&\hpi = \argmax_{\pi \in \Pi} \cb{\frac{1}{n} \sum_{i = 1}^n \lambda_i H_i (2\pi(X_i) - 1)}, \ \ \lambda_i = \abs{\hGamma_i}, \ \ H_i =  \text{sign}\p{\hGamma_i}, 
\end{split}
\end{equation}
where we train a classifier $\pi(\cdot)$ with response $H_i$ using sample weights $\lambda_i$.
Given this formalism, we can build on existing tools for
weighted classification to learn $\hpi$; see \citet*{zhou2018offline} for a further
discussion.\footnote{Some popular approaches for solving problems of the form \eqref{eq:classif}
include best-subset empirical risk minimization \citep{chen2016best} and
optimal trees \citep{bertsimas2017optimal}. Due to the computational difficulty
of solving the problem \eqref{eq:classif} exactly, it may also be of interest to consider the empirical performance of 
alternative methods that solve an approximation to our weighted classification
problem, e.g., support vector machines \citep{cortes1995support} or
recursive partitioning \citep*{breiman1984classification}.
However, we caution that our formal results only apply to methods that solve the problem
\eqref{eq:classif} exactly; see \citet{wager2020regression} for further discussion.}
In all our experiments, we set $\Pi$ to be a class of finite-depth decision trees (see Section \ref{sec:donsker}
for a definition), and solve the optimization problem in Step 3 using our companion \texttt{R}-package \texttt{policytree}
\citep{sverdrup2020policytree,CRAN}; see \citet{zhou2018offline} for further details and motivation behind the computational
strategy taken in this package.

\subsection{The California GAIN Program}
\label{sec:gain}

The Greater Avenues for Independence (GAIN) program, started in 1986, is a welfare-to-work program that provides participants with a mix of educational
resources and job search assistance. Between 1988 and 1993, the Manpower Development Research Corporation
conducted a randomized study to evaluate the program. As described in \citet*{hotz2006evaluating},
randomly chosen registrants were eligible to receive GAIN benefits immediately, whereas others
were embargoed from the program until 1993.
All experimental subjects were followed for a 9-year post-randomization period and,
as documented by \citet{hotz2006evaluating}, eligibility for GAIN had a significant impact
on mean quarterly income averaged over this 9-year period.

Our current question is
whether we can find ways to prioritize treatment to some subgroups of GAIN registrants particularly
likely to benefit from it. We consider data from four counties, Alameda, Riverside, Los Angeles
and San Diego, resulting in $n = 19,170$ observations, and use $p = 28$ covariates,
including demographics, education, and per-quarter earnings for 10 quarters preceding treatment.
As in \citet{hotz2006evaluating}, we use average quarterly income over the 9-year post-randomization period (in \$1000s)
as our outcome.

Each county participating in the GAIN evaluation conducted its own randomized controlled trial,
and the counties had considerable freedom in how they carried out the randomization. In particular,
counties had flexibility in choosing whom to enroll in the randomized trial, and which
fraction of participants to randomize into treatment. The data reflects this heterogeneity in study
specifications: The per-county average outcome for controls varied from 0.64 to 1.04 thousand
dollars per quarter, while the per-county fraction of treated units varied from 0.50 to 0.86.

We use this dataset to design a semi-synthetic observational study by pooling the data
from all four counties under consideration. Because the mean control outcome and treatment
fraction vary from county to county (and are in fact correlated), we expect that an uncorrected
analysis of the pooled data would suffer from confounding. In an attempt to correct for the
confounding that arises from pooling we pursue a selection-on-observables
strategy, and assume that controlling for the $p = 28$ covariates described above is enough to
correct for the different study specifications used in different counties.

Our method starts by computing doubly robust scores for the treatment effect,
and learning policies by empirical maximization as in \eqref{eq:meta_estimator}.
We use the augmented inverse-propensity weighted scores of \citet*{robins1994estimation},
with nuisance component estimates from
generalized random forests \citep*{athey2018generalized,breiman2001random},\footnote{The one major deviation between how
we compute scores below and the assumptions of Theorem \ref{theo:main} is that, here, we use leave-one-out
(or out-of-bag) estimates for $\tau(X_i)$, etc., whereas Theorem \ref{theo:main} assumed $K$-fold estimation.
The reason for this choice is that, as discussed in \citet{breiman2001random}, random forests are particularly
well suited for leave-one-out estimation, and allow the analyst to obtain such estimates at essentially no additional
computational cost.}
\begin{align}
\label{eq:grf-opt}
&\hpi = \argmax\cb{\frac{1}{n} \sum_{i = 1}^n \p{2\pi(X_i) - 1} \p{\hGamma_i - C} : \pi \in \Pi}, \\
\label{eq:grf-score}
&\hGamma_i = \htau^{(-i)}(X_i) \\
\notag
& \ \ \ \ + \frac{W_i - \he^{(-i)}(X_i)}{\he^{(-i)}(X_i) \p{1 - \he^{(-i)}(X_i)}} \p{Y_i - \hat{f}^{(-i)}(X_i) - (W_i - \he^{(-i)}(X_i))\htau^{(-i)}(X_i)},
\end{align}
where \smash{$\hat{f}(x)$} and \smash{$\he(x)$} are random forest estimates of $\mathbb{E}[Y_i \cond X_i= x]$ and
$\mathbb{E}[W_i \cond X_i= x]$ respectively, \smash{$\htau(\cdot)$} is an causal forest\footnote{Random
forests are a type of adaptive nearest neighbor estimator that use an ensemble of trees to define a relevant
neighborhood function for each query point; see \citet{athey2018generalized} for a discussion.
Causal forests use the adaptive neighborhood function implied by a forest to fit a partially linear model using the method of
\citet{robinson1988root}; see \citet{nie2017learning} for formal results motivating the use of local partially linear
modeling for heterogeneous treatment effect estimation, and Section 1.3 of \citet{athey2019estimating} for a
discussion of how this partially linear modeling is carried out in causal forests.
We emphasize that, for our purposes, random forests are simply used as a convenient non-parametric estimator of relevant
nuisance components, specifically $f(x)$ and $e(x)$ here, and could
seamlessly be replaced with other methods such as boosting or neural networks. The shape of the learned policy
$\hpi$ is determined in the optimization step 3, which only depends on the random forests through the predictions
used to form doubly robust scores \smash{$\hGamma_i$}.}
estimate of the conditional average treatment effect, and $C$ is a parameter measuring the cost of
treatment. Tuning parameters for all forests were selected by leave-one-out cross-validation.\footnote{The regression surfaces
\smash{$\hat{f}(x)$} and \smash{$\he(x)$} were tuned to optimize mean-squared error.
As advocated in \citet{nie2017learning}, the conditional average treatment effect function was tuned to
optimize the error of a local residual-on-residual regression.}
Here, we set $C = 0.14$ to roughly match the average treatment effect with the goal of ensuring that the optimal treatment
rule is not trivial (i.e., we can only achieve non-zero utility gains by exploiting treatment heterogeneity).

Before starting to optimize policies we first run a brief sanity check on our selection-on-observables strategy,
and confirm the ability of estimators that build on this assumption to accurately recover the average treatment
effect we would get using a proper randomization-based estimator that does not pool data across counties.
The natural doubly robust estimator of the average treatment
effect in our setting is \smash{$\htheta_{DR} = \sum_{i = 1}^n \hGamma_i / n$}, with scores
\smash{$\hGamma_i$} as in \eqref{eq:grf-score}. We compare it to a naive difference-in-means estimator
$\htheta_{DM} = \text{avg}\cb{Y_i : W_i = 1} - \text{avg}\cb{Y_i : W_i = 0}$ that does not attempt to correct for
bias due to pooling, and to an ``oracle'' doubly robust estimator that
does not estimate propensity scores from covariates but instead uses the true per-county treated fractions:
\smash{$\htheta_{DR}^* = \sum_{i = 1}^n \hGamma_i^* / n$} with
\begin{equation}
\label{eq:grf-score-oracle}
\begin{split}
&\hGamma_i^* = \htau^{(-i)}(X_i) + \frac{W_i - \he^*_i}{\he^*_i \p{1 - \he^*_i}} \p{Y_i - \hat{f}^{(-i)}(X_i) - (W_i - \he^*_i)\htau^{(-i)}(X_i)}, \\
&\he_i^* = \sum_{j = 1}^n W_j 1\p{\cb{G_j =  G_i}} \, \Big/\,\sum_{j = 1}^n 1\p{\cb{G_j =  G_i}},
\end{split}
\end{equation}
where $G_i \in \cb{\text{Alameda}, \, \text{Riverside}, \, \text{Los Angeles}, \, \text{San Diego}}$ denotes
the county-membership of the $i$-th sample. Because \smash{$\htheta_{DR}^*$} uses the true per-county treatment fractions
$\he_i^*$ and estimates nuisance components using cross-fitting, the point estimates will be $\sqrt{n}$-consistent
and the associated confidence intervals asymptotically valid essentially without assumptions \citep*{rothe2018flexible,wager2016high}.
The resulting point estimates for the average treatment effect ($\pm 1$ standard error)
are: $\htheta_{DR} = 0.141 \pm 0.026$ for the feasible doubly-robust estimator, $\htheta_{DR}^* = 0.146 \pm 0.028$ for
the oracle doubly-robust estimator, and  $\htheta_{DM} = 0.208 \pm 0.028$ for the naive difference in means.
Thus, it appears that pooling county information results in confounding,
but that controlling for available covariates helps.

\begin{table}[t]
\begin{center}
\begin{tabular}{r|cc|}
 & non-white & white \\ \hline
fraction treated & 76\% & 81\%  \\
mean control outcome & 0.79 & 0.90 \\ \hline
\end{tabular}
\caption{Outcome is mean quarterly income (in \$1000) averaged over 9 years post-intervention.
Differences in mean responses between white and non-white respondents are both significant at the
$0.05$ level using a Welch two-sample $t$-test.}
\label{tab:conf}
\end{center}
\end{table}

We now move to learning a policy $\hpi$. In doing so, however, we note that
caution is warranted because we have measured features pertaining to race, ethnicity, age and gender.
On the one hand, there may be legal restrictions on the use of these features for treatment allocation but, on the other hand, they appear to act as
counfounders. For example, as shown in Table \ref{tab:conf}, white GAIN registrants were randomized to treatment
at higher rates than non-white registrants, and also white controls had higher outcomes than non-white controls.
Our approach allows us to seamlessly use such sensitive variables for deconfounding without using them for
policy allocation: We use these variables when estimating the nuisance components in \eqref{eq:grf-score}, but then
omit them from the maximization step \eqref{eq:grf-opt} that produces the policy.

\begin{figure}[t]
\begin{center}
\begin{tabular}{cc}
\underline{\it depth 1 policy} & \underline{\it depth 2 policy} \\[3.5mm]
\begin{tikzpicture}[
   level distance=2cm,sibling distance=2cm,
   edge from parent path={(\tikzparentnode) -- (\tikzchildnode)},
   every node/.append style={align=center}]
\Tree[.{is high school graduate} 
  \edge node[auto=left] {\it no};
         [ .{\textit{don't treat}  \\[1.5mm] ($n = 9,477$)} ]
  \edge node[auto=right] {\it yes};
		[ .{\textit{treat}  \\[1.5mm] ($n = 9,693$)}  ]]
\end{tikzpicture} &
\begin{tikzpicture}[
   level distance=1.5cm,
   edge from parent path={(\tikzparentnode) -- (\tikzchildnode)},
   every node/.append style={align=center}]
\Tree[.{was paid 3 quarters ago}
    \edge node[auto=left, pos=0.35] {\it no};
          [.{has children} [ .{\textit{don't treat} \\[1.5mm]  ($n = 2,486$)} ]
               [ .{\textit{treat} \\[1.5mm]  ($n = 2,757$)}  ]]
    \edge node[auto=right, pos=0.35] {\it yes};           
          [.{is high school graduate} [ .{\textit{don't treat} \\[1.5mm] ($n = 7,325$)}  ]
               [ .{\textit{treat} \\[1.5mm]  ($n = 6,602$)}  ]]]
\end{tikzpicture}
\end{tabular}
\caption{Example of optimal depth-1 and -2 policy trees learned by optimizing the augmented
inverse-propensity weighting loss function.}
\label{fig:tree}
\end{center}
\end{figure}
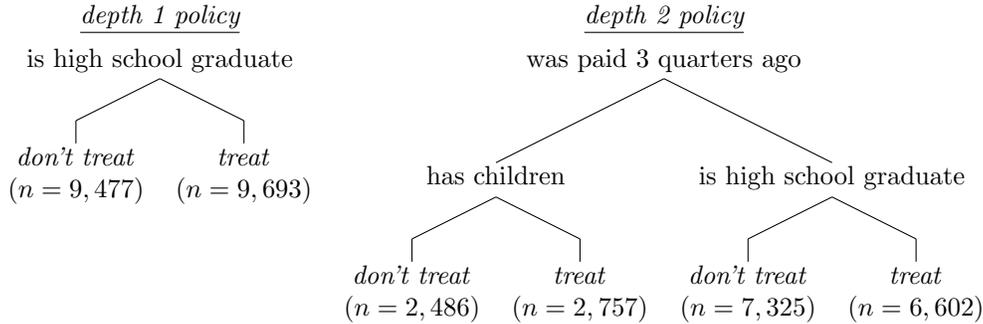

For our policy class $\Pi$, we consider decision trees of depth either 1 or 2.
The learned decision rules are shown in Figure \ref{fig:tree}. Interestingly, the depth-1 and 2
trees make the same decisions for the roughly 3/4 of GAIN registrants who were paid 3 quarters prior to
randomization, but the depth-2 tree chooses to switch to a different rule for those who weren't paid 3 quarters prior.

\begin{table}[t]
\centering
\begin{tabular}{r|cc|}
& \multicolumn{2}{c|}{estimated improvement} \\
method & fitted propensities & true propensities \\
  \hline
  plug-in & $ 0.077 \pm 0.026 $ & $ 0.063 \pm 0.028 $ \\ 
  IPW depth 1 & $ 0.065 \pm 0.026 $ & $ 0.048 \pm 0.028 $ \\ 
  IPW depth 2 & $ 0.043 \pm 0.026 $ & $ 0.029 \pm 0.028 $ \\ 
  AIPW depth 1 & $ 0.068 \pm 0.026 $ & $ 0.050 \pm 0.028 $ \\ 
  AIPW depth 2 & $ 0.091 \pm 0.026 $ & $ 0.080 \pm 0.028 $ \\ 
   \hline
\end{tabular}
\caption{Estimate of the utility improvement of various policies over a random assignment baseline,
$\pm 1$ standard error.
The plug-in policy simply thresholds causal forest predictions at \smash{$\htau^{(-i)}(X_i) > C$},
the inverse-propensity weighted trees (IPW) are following \citet{kitagawa2015should}, and
the trees scored via augmented inverse propensity-weighting (AIPW) are instances of the method studied here.
The left column estimates improvement via \eqref{eq:CV_advantage}, whereas the right column brings in
county membership information to obtain a randomization-based estimator of improvement \eqref{eq:CV_advantage_oracle}.}
\label{tab:holdout_advantage}
\end{table}

In order to choose tree depth and, more broadly, to evaluate the accuracy of the policy learning procedure, we
recommend cross-validation. We randomly divide the data into $K$ folds $\set_k$, $k = 1, \, ..., \, K$ and, for each
fold, learn a policy $\hpi^{(-k)}(\cdot)$ using all but the data in $\set_k$. Here, we use $K = 10$.
Finally, we estimate improvement over a random baseline as
\begin{equation}
\label{eq:CV_advantage}
\hA_{CV} = \frac{1}{n} \sum_{k = 1}^K \sum_{i \in \set_k} \p{2 \hpi^{(-k)}(X_i) - 1} \, \hGamma_i.
\end{equation}
Table \ref{tab:holdout_advantage} shows the estimated improvement of our depth-1 and -2 trees, as well as two
baselines: A variant of the inverse-propensity weighted method of \citet{kitagawa2015should} using the propensities
used to construct \eqref{eq:grf-score}, as well as a plug-in policy that does not obey our functional form restriction,
and simply treats all samples with \smash{$\htau^{(-i)}(X_i) > C$}. Our depth-2 trees achieve markedly better performance
than the depth-1 trees. Interestingly, the depth-2 tree is also competitive with the unconstrained plug-in
estimator. Based on this analysis, we prefer the depth-2 tree in Figure \ref{fig:tree}.

One potential concern with this analysis is that our evaluation hinges on validity of the selection-on-observables assumption,
as well as accuracy of the doubly robust scores \smash{$\hGamma_i$} from \eqref{eq:grf-score}. To assuage this concern,
we also computed a version of the improvement measure \eqref{eq:CV_advantage}, but with scores \smash{$\hGamma_i^*$}
computed using the true per-county treatment fractions as in \eqref{eq:grf-score-oracle}:
\begin{equation}
\label{eq:CV_advantage_oracle}
\hA^*_{CV} = \frac{1}{n} \sum_{k = 1}^K \sum_{i \in \set_k} \p{2 \hpi^{(-k)}(X_i) - 1} \, \hGamma_i^*.
\end{equation}
As seen in the second rightmost column of Table \ref{tab:holdout_advantage}, our feasible evaluation discussed above gave
the correct ordering for the methods, but was somewhat optimistic in terms of the quality of the learned policies.
The formal properties of treatment rules whose complexity is tuned via cross-validation are considered by
\citet{mbakop2016model}.\footnote{Recall that cross-validation is a means of evaluating the quality of the policy learning
procedure, not the decision that was produced by a specific realization of the procedure. If we want an accuracy assessment
that is valid conditionally on the learned rule \smash{$\hpi(\cdot)$}, one can either use a single test-train split, or use the more
sophisticated data carving approach of \citet*{fithian2014optimal}.}

\subsection{Simulation Study with Binary, Endogenous Treatments}
\label{sec:simu}

In order to develop a richer quantitative understanding of the behavior of our method, we
now turn to a simulation study. Here, we consider a setting with a binary, endogenous treatment
$W_i$ and a binary instrument $Z_i$ and assume homogeneity as in \eqref{eq:iv_homog}. In this
case, our method chooses the policy
\smash{$\hpi = \argmax\{\frac{1}{n} \sum_{i = 1}^n \p{2\pi(X_i) - 1} \hGamma_i : \pi \in \Pi\}$},
where \smash{$\hGamma_i$} is a cross-fit doubly robust score with estimates of the compliance weights
as in \eqref{eq:iv_weight}:
\begin{equation}
\label{eq:iv_score}
\begin{split}
&\hGamma_i = \htau^{(-i)}(X_i) + \hg^{(-i)}(X_i, \, Z_i) \p{Y_i - \hat{f}^{(-i)}(X_i) - (W_i - \he^{(-i)}(X_i))\htau^{(-i)}(X_i)}, \\
&\hg^{(-i)}(X_i, \, Z_i) = \frac{1}{\hDelta^{(-i)}(X_i)} \frac{Z_i - \hz^{(-i)}(X_i)}{\hz^{(-i)}(X_i)(1 - \hz^{(-i)}(X_i))} , \\
\end{split}
\end{equation}
\sloppy{where \smash{$\Delta(x) = \PP{W_i = 1 \cond Z_i = 1, \, X_i = x} -  \PP{W_i = 1 \cond Z_i = 0, \, X_i = x}$} is the conditional average
effect of the instrument on the treatment, \smash{$z(x) = \PP{Z_i = 1 \cond X_i = x}$}, \smash{$f(x) = \EE{Y_i \cond X_i = x}$},
\smash{$e(x) = \PP{W_i = 1 \cond X_i = x}$}, and $\tau(x)$ is the conditional average treatment effect as specified
in \eqref{eq:iv_homog}. We estimate all nuisance components via
random forest methods with the package \texttt{grf}, and use an instrumental forest for $\tau(\cdot)$, a
causal forest for $\Delta(\cdot)$, and a regression forest for $f(\cdot)$, $e(\cdot)$ and $z(\cdot)$.}

\begin{figure}
\begin{center}
\begin{tabular}{ccc}
\includegraphics[width=0.45\textwidth, trim=8mm 8mm 8mm 8mm]{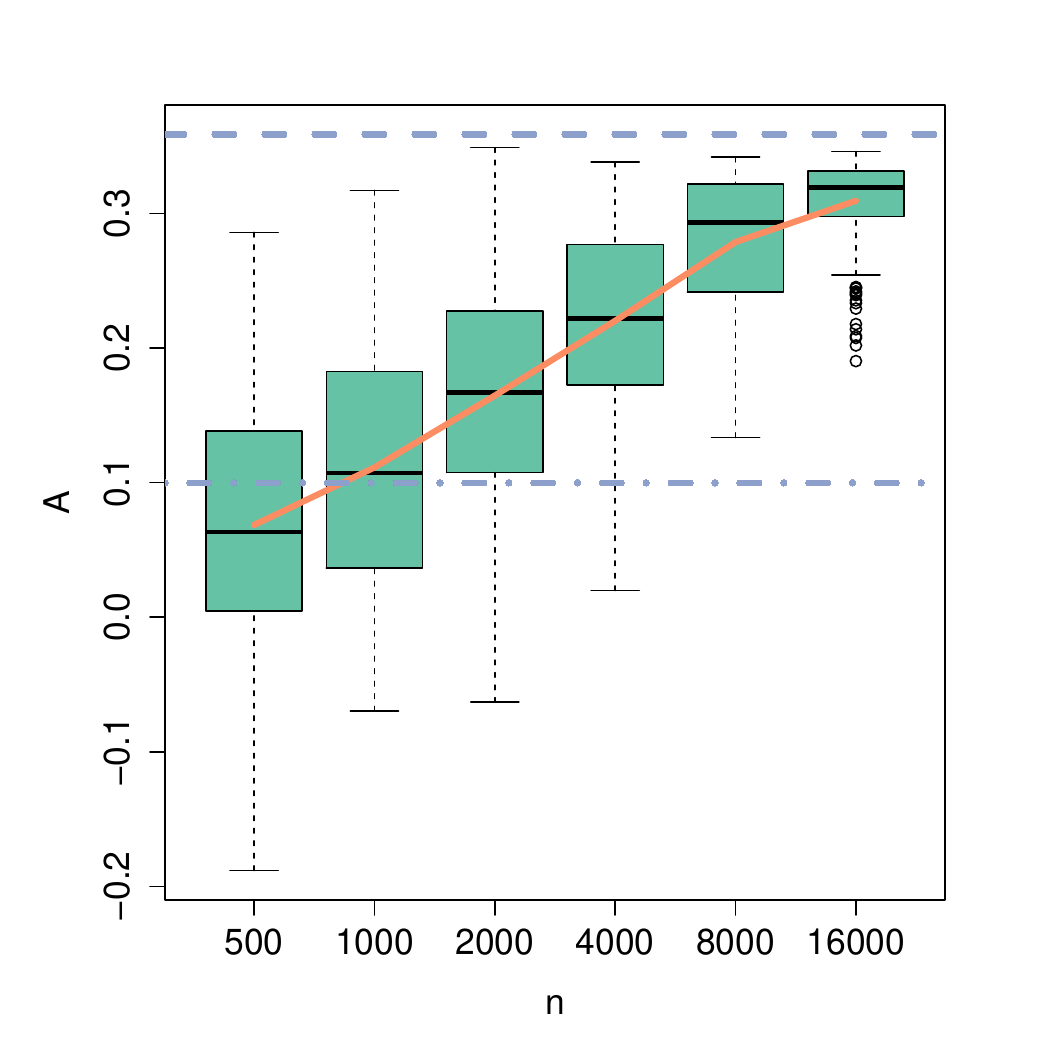} &
\includegraphics[width=0.45\textwidth, trim=8mm 8mm 8mm 8mm]{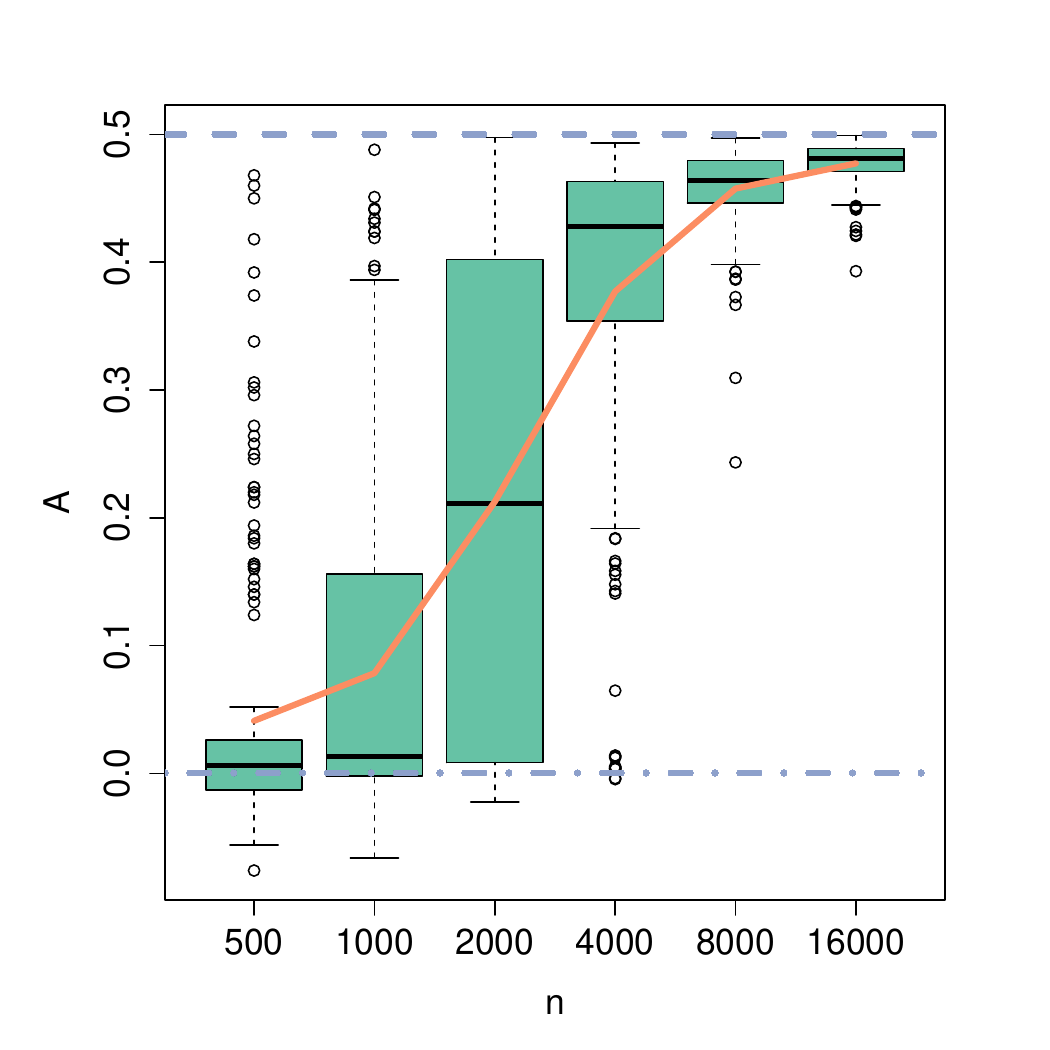} \\[2mm]
$\tau(\cdot)$ as in \eqref{eq:tau_add} & $\tau(\cdot)$ as in \eqref{eq:tau_prod}
\end{tabular}
\caption{Distribution of the improvement $A(\pi) = \EE{(2\pi(X_i) - 1)\tau(X_i)}$ for policies learned
by optimizing the scores \eqref{eq:iv_score} over the class $\Pi$ of depth-2 trees, for different values
of sample size $n$. Each box plot summarizes the distribution of \smash{$A(\hpi)$} over 200
simulation replications, while the solid line shows the average of \smash{$A(\hpi)$}. The lower
horizontal line shows $A(\pi)$ for the best policy that does not use the features $X_i$ (i.e., either
always treat or never treat), and the upper horizontal line shows the supremum of $A(\pi)$ over
the class $\Pi$.}
\label{fig:instrumental_simu}
\end{center}
\end{figure}

In this simulation experiment, we generate data independently as follows, for various
choices of $n$ and $\tau(\cdot)$:
\begin{equation}
\begin{split}
&X \sim \nn\p{0, \, \ii_{10 \times 10}}, \ \ Z \cond X \sim \text{Bernoulli}\p{1/\p{1 + e^{-X_3}}}, \\
&\varepsilon \cond X, \, Z \sim \nn\p{0, \, 1}, \ \ Q \cond X, \, Z, \, \varepsilon \sim \text{Bernoulli}\p{1/\p{1 + e^{-\varepsilon - X_4}}}, \\
&W = Q \land Z, \ \ Y = \p{X_3 + X_4}_+ + W \tau(X) + \varepsilon.
\end{split}
\end{equation}
Note that $W$ is in fact endogenous, because $Q$ (and thus also $W$) is more likely to be 1
when the noise term $\varepsilon$ is large. Given this setup, we consider $\tau(\cdot)$ functions
\begin{align}
\label{eq:tau_add}
&\tau(x) = \p{\p{x_1}_+ + \p{x_2}_+ - 1} / 2 \ \eqand\\
\label{eq:tau_prod}
&\tau(x) = \text{sign}\p{x_1 x_2} / 2.
\end{align}
In both cases, we learn $\pi(\cdot)$ over the class $\Pi$ of depth-2 trees and note that
best non-parametric policy \smash{$\pi^*(x) = 1\p{\cb{\tau(x) > 0}}$} belongs to $\Pi$ in
case \eqref{eq:tau_prod} but not in case \eqref{eq:tau_add}.

In Figure \ref{fig:instrumental_simu}, we display the improvement $A(\pi) = \EE{(2\pi(X_i) - 1)\tau(X_i)}$ of our
learned policies relative to a random assignment baseline, for different values of $n$. Over all, we see that
the regret of the learned policies improves with $n$, and approaches best-in-class regret as $n$ gets large.
We also note an interesting difference in the behavior of the learned rules in settings \eqref{eq:tau_add} and
\eqref{eq:tau_prod}. In the first case, $\tau(\cdot)$ is continuous, and regret improves smoothly with sample size.
Conversely, in the second case where $\tau(\cdot)$ has sharp jumps, we observe something of a phase transition
between $n = 2,000$ and $n = 4,000$, as our trees become able to consistently make splits that roughly match
the jumps in $\tau(\cdot)$.

\section{Discussion}

In this paper, we proposed an approach to policy learning in the observational study setting that builds on classical
ideas for semiparametrically efficient treatment effect estimation. Our main result is that doubly robust estimators
of average treatment effects can be adapted for policy evaluation, and that the policy that maximizes the resulting
doubly robust value estimate over a pre-specified class $\Pi$ satisfies rate-optimal guarantees for minimax regret.
Our approach decouples estimation of nuisance components used for the doubly robust scores from optimization
of the doubly robust value function, and thus allows practitioners flexibility in how they implement each step.

Our formal discussion focused on regret bounds for policy learning. A natural follow-up question is to ask
for confidence sets guaranteed to contain an optimal policy: For example, if $\Pi$ is the set of depth-$L$
decision trees, can we identify a subset of $\Pi$ guaranteed to contain a value-maximizing policy in $\Pi$ 
with high probability? Some early results in this direction are reported by \citet{rai2018statistical}. Meanwhile,
\citet{armstrong2015inference} consider the related task of identifying a subset of the population
we are confident will benefit from the policy intervention.

Another natural direction to extend our results is towards dynamic decision making problems, where
the policy maker needs to make a sequence of decisions, potentially depending on time-varying
covariates. The problem of doubly robust policy evaluation in this setting has been considered by
\citet*{thomas2016data} and \citet*{zhang2013robust}. \citet*{nie2019learning} proposed
a method for learning observational stopping rules from observational data that is both computationally
feasible and robust to confounding. Obtaining a more comprehensive landscape of the problem of dynamic
policy learning in observational studies would be of considerable interest.

Finally, all results presented here relied on point-identification of treatment effects, either via a selection on
observables assumption or via an instrument that satisfies conditional homogeneity. Some applications, however,
do not allow for such clean assumptions, and thus call for methods for policy learning that are robust to failures
of identifying assumptions. \citet{kallus2018confounding} consider the problem of policy learning under an
approximate selection-on-observables assumption in the sense of \citet{rosenbaum2002observational}.
It would also be of interest to study what can be done if we only have access to a monotone instrument,
as in \citet{manski2000monotone}.



\newpage

\begin{appendix}

\ifecma
\title{Supplemental Material: \\ Policy Learning with Observational Data}
\else
\title{Supplemental Material}
\author{}
\fi

\date{}
\maketitle

\ifecma
\else
\vspace{-5\baselineskip}
\fi

\section{Characterizing the VC Dimension}
\label{app:VC}

As a preliminary to our technical argument, we start by reviewing some practical characterizations
of the VC dimension in terms of covering numbers in Hamming distance.
For any discrete set of points $\cb{X_1, \, ..., \, X_m}$ and any $\varepsilon > 0$, define the
$\varepsilon$-Hamming covering number $N_H(\varepsilon, \, \Pi, \, \cb{X_1, \, ..., \, X_m})$
as the smallest number of policies $\pi : \cb{X_1, \, ..., \, X_m} \rightarrow \cb{0, \, 1}$
(not necessarily contained in $\Pi$)
required to $\varepsilon$-cover $\Pi$ under Hamming distance,
\begin{equation}
H(\pi_1, \, \pi_2) = \frac{1}{m} \sum_{j = 1}^m \ind\p{\cb{\pi_1(X_j) \neq \pi_2(X_j)}}. 
\end{equation}
Then, define the $\varepsilon$-Hamming entropy of $\Pi$ as $\log\p{N_H(\varepsilon, \, \Pi)}$, where
\begin{equation}
N_H(\varepsilon, \, \Pi) = \sup\cb{N_H(\varepsilon, \, \Pi, \, \cb{X_1, \, ..., \, X_m}): \, X_1, \, ..., \, X_m \in \xx; \,  m \geq 1} 
\end{equation}
is the number of functions needed to $\varepsilon$-cover $\Pi$ under
Hamming distance for any discrete set of points.
We note that this notion of entropy is purely geometric, and does not depend on the
distribution used to generate the $X_i$.

As argued in \citet{pakes1989simulation}, a class $\Pi$ has a finite VC dimension if and
only if there is a constant $\kappa$ for which
\begin{equation}
\label{eq:entropy_decay}
\log\p{N_H(\varepsilon, \, \Pi_n)} \leq \kappa \log\p{\varepsilon^{-1}}
\ \text{ for all } \
0< \varepsilon < \frac{1}{2}.
\end{equation}
Moreover, there are simple quantitative bounds for Hamming entropy in terms of the VC dimension:
If $\Pi$ is a VC class of dimension $\VC(\Pi)$, then \citep{haussler1995sphere} 
\begin{equation}
\label{eq:haussler}
\begin{split}
\log\p{N_H(\varepsilon, \, \Pi)} 
&\leq \VC(\Pi)\p{\log\p{\varepsilon^{-1}} + \log(2) + 1} + \log\p{\VC(\Pi) + 1} + 1 \\
&\leq 5 \VC(\Pi)\log\p{\varepsilon^{-1}}\ \text{ for all } \
0< \varepsilon < \frac{1}{2}
\end{split}
\end{equation}
whenever $\VC(\Pi) \geq 2$. 
Conversely, recall that if $\Pi$ has VC-dimension $d$ it can shatter a set
of $d$ points, and so we must have $N_H(1/d, \, \Pi) \geq 2^d$. Thus, the VC dimension $d$ of any
class whose Hamming entropy satisfies \eqref{eq:entropy_decay} must be bounded via the relationship
\begin{equation}
d \log(2) \leq \kappa \log(d).
\end{equation}
Whenever we invoke Assumption \ref{assu:VC} in our proof, we actually work in terms of the
covering number bound \eqref{eq:haussler} and assume that $\VC(\Pi) \geq 2$ (the case with
$\VC(\pi) = 1$, corresponding to non-personalized decision rules, is trivial).

\section{Additional Simulation Experiments}
\label{app:simu_continuous}

We complement our experiments from Section \ref{sec:experiments} with another
simulation example where, now, the treatment dose $W_i \in \RR$ is continuous.
As discussed in Section \ref{sec:identification}, we consider policies that infinitesimally nudge the treatment
dose $W_i$ for select samples; the value $V(\pi)$ of a policy $\pi$ is then:
\begin{equation}
\pi : \xx \rightarrow \cb{0, \, 1}, \ \ \ \ \ V(\pi) = \EE{\pi(X_i) \p{ \sqb{\frac{d}{d\nu} Y_i(W_i + \nu)}_{\nu = 0} - C}}, 
\end{equation}
where $C$ is a cost of treatment.
We assume $W_i$ to be exogenous. As always, we learn our policy $\hpi$ via
\smash{$\hpi = \argmax\{\frac{1}{n} \sum_{i = 1}^n \p{2\pi(X_i) - 1} (\hGamma_i - C) : \pi \in \Pi\}$},
and the \smash{$\hGamma_i$} are appropriate cross-fit doubly robust scores \eqref{eq:continuous_riesz},
\begin{equation}
\label{eq:Gamma_cont}
\begin{split}
&\hGamma_i = \sqb{\frac{d}{dw} \hatm^{(-i)}\p{X_i, \, w}}_{w = W_i} \\
&\ \ \ \ \ \ \ \ \ \ - \frac{d}{dw} \sqb{\log\p{\hf^{(-i)}\p{w \cond X_i}}}_{w = W_i} \p{Y_i -  \hatm^{(-i)}\p{X_i, \, W_i}},
\end{split}
\end{equation}
where $f(\cdot \cond x)$ denotes the conditional density of $W_i$ given $X_i = x$,
and $m(x, \, w) = \EE{Y_i \cond X_i = x, \, W_i = w}$.

Unlike in our previous examples, the non-parametric regression problems underlying \eqref{eq:Gamma_cont}
have not received much attention in the statistical learning literature. First, \eqref{eq:Gamma_cont} requires
estimating derivatives of conditional response-functions; but many popular machine learning methods, such as
random forests or boosted trees, do not have differentiable predictive surfaces. Second, the problem of estimating
a conditional density function $f(\cdot \cond x)$ presents its own numerical challenges.

Here, we approach the problem as follows. In order to make sure that the derivatives of $\hatm(\cdot)$ and $\hf(\cdot)$ are
good estimates of $m(\cdot)$ and $f(\cdot)$ respectively, we use penalized series estimators throughout.
We fit $\hatm(X_i, \, W_i)$ by penalized regression on 3rd-order Hermite polynomials in $(X_i, \, W_i)$.
Meanwhile, we fit the conditional density function $f(\cdot \cond X_i)$ by adapting Lindsey's method,
a technique for estimating distribution functions using software for generalized linear modeling
\citep{efron1996using,lindsey1974comparison}. In the case without covariates,
Lindsey's method involves first discretizing the support of $W_i$ into a union of non-overlapping equal-length intervals and,
as with a histogram, counting the number of samples $W_i$ that fall within each interval. Then, these histogram
counts are fit via Poisson regression using a series expansion of $W_i$. As shown in \citet{efron2011tweedie}, the
log-derivative of the estimated density function is well-behaved as an estimate of the log-derivative of the
true density. Now, in the case with covariates, we again discretize the support of $W_i$ into $K$ non-overlapping intervals.
However, instead of making a histogram, we now duplicate each sample $K$ times: For each sample $i = 1, \, ..., \, n$ and interval $k = 1, \ ..., \, K$
we create a datapoint $(X_i, \, w_k, \, L_{ik})$, where $w_k$ is the mid-point of the $k$-th interval and $L_{ik}$ is an
indicator for whether the $W_i$ is in the $k$-th interval. Finally, we fit this model by penalized logistic regression
on full interactions between 3rd-order Hermite polynomials in $X_i$ and an appropriate basis expansion $b(w)$
in $w$ discussed further below. In all cases, we fit penalized regression via \texttt{glmnet} \citep*{friedman2010regularization},
with the amount of penalization tuned via cross-validation.

We consider the following simulation designs, loosely motivated by a probit choice model in a pricing application
(i.e., where $W_i$ acts as a price and $Y_i$ is a choice to purchase). In all cases, we generate independent samples as
below, with $p = 6$:
\begin{equation}
\label{eq:cont_setting}
\begin{split}
&X_i \sim \nn\p{0, \, \ii_{p \times p}}, \ \ \ \ U_i = 5 \,\Big/ \p{1 + 3 e^{-(X_{i1} + X_{i2})}} - 0.5,   \\
&W_i \cond X_i \sim \law_w(X_i), \ \ \ \  Y_i \cond U_i, \, W_i \sim \text{Bernoulli}\p{\Phi(W_i - U_i)},
\end{split}
\end{equation}
where $\Phi(\cdot)$ is the standard Gaussian cumulative distribution function. We consider two choices for the
conditional distribution $\law_w$ of $W_i$ conditionally on $X_i$:
\begin{align}
\label{eq:cont_W1}
&\text{Gaussian:} &W_i = 3 \,\Big/ \p{1 + 3e ^{-(X_{i1} + X_{i3})}} + \varepsilon_i, \ \ \ \ \varepsilon_i \cond X_i \sim \nn\p{0, \, 1}, \ \eqand \\
\label{eq:cont_W2}
&\text{Non-Gaussian:} &W_i = 3 \,\Big/ \p{1 + 3e ^{-(X_{i1} + X_{i3} + \eta_i)}} + \varepsilon_i, \ \ \ \ (\varepsilon_i, \, \eta_i) \cond X_i \sim \nn\p{0, \, \ii_{2 \times 2}}.
\end{align}
In principle, the Gaussian case appears substantially easier than the non-Gaussian case, because the logistic regression
problem underlying Lindsey's method as above is well-specified with a quadratic expansion in $w$, i.e., $b(w) = (1 \ \ w \ \ w^2)$. 
In the non-Gaussian case, no similar simplifications apply. In our experiments, we in fact set $b(w)$ to be the quadratic
expansion in the Gaussian case; in the non-Gaussian case, we set $b(w)$ to a 5th order natural spline basis.

Before evaluating the accuracy of policy learning in this setting, we present some performance diagnostics on the associated
doubly robust average derivative estimator $\htheta_{DR} = \sum_{i = 1}^n \hGamma_i / n$ as, despite attracting a
fair amount of interest in the literature on asymptotic estimation \citep*[including][]{chernozhukov2016locally,chernozhukov2018double,hirshberg2018balancing}, we are not aware of
existing Monte Carlo evaluations of this estimator in the literature.\footnote{The closest experiments we are aware from
are from \citet{graham2018semiparametrically} and \citet{hirshberg2018balancing}, who report results results for doubly
robust estimation in a closely related (but more restricted) model with a conditionally linear specification 
$\EE{Y_i \cond X_i = x, \, W_i = w} = m(x) + w \tau(x)$.} We report bias and root-mean squared error for the doubly robust
estimator $\htheta_{DR}$, the pure regression estimator $\htheta_{reg} = \sum_{i = 1}^n d/dw \, \hatm^{(-i)}(X_i, \, W_i) / n$,
and the pure weighting estimator $\htheta_{weight} = \sum_{i = 1}^n d/dw \log \hf^{(-i)}(X_i, \, W_i) \, Y_i / n$. We also
report mean-squared standardized error $S = \mathbb{E}[(\htheta_{DR} - \theta)^2 / \hsigma^2]^{1/2}$ with
$\hsigma^2 = \sum_{i = 1}^n \hGamma_i / (n (n-1))$
which, under the conditions of Assumption \ref{assu:LR}, should converge as $\limn S = 1$.

\begin{table}[t]
\centering
\begin{tabular}{|r|c|cc|cc|ccc|c|}
 \hline
   &    & \multicolumn{2}{c|}{regression} & \multicolumn{2}{c|}{weighted} & \multicolumn{3}{c|}{doubly robust} & policy \\
  & $n$ & bias & RMSE &  bias & RMSE &  bias & RMSE & $S$ & value \\ 
    \hline
 \multirow{4}{*}{\rotatebox[origin=c]{90}{setup 1}} 
   & 600 & -0.056 & 0.058 & -0.132 & 0.133 & -0.035 & 0.037 & 4.59 & 0.014 \\ 
   & 1800 & -0.035 & 0.036 & -0.095 & 0.096 & -0.017 & 0.019 & 2.97 & 0.024 \\ 
   & 5400 & -0.022 & 0.022 & -0.081 & 0.081 & -0.010 & 0.010 & 2.60 & 0.028 \\ 
   & 16200 & -0.012 & 0.013 & -0.073 & 0.073 & -0.006 & 0.006 & 2.54 & 0.029 \\ 
   \hline
  \multirow{4}{*}{\rotatebox[origin=c]{90}{setup 2}}
    & 600 & -0.069 & 0.072 & -0.062 & 0.063 & -0.049 & 0.050 & 8.01 & 0.018 \\ 
    & 1800 & -0.040 & 0.041 & -0.052 & 0.053 & -0.026 & 0.027 & 6.26 & 0.033 \\ 
    & 5400 & -0.023 & 0.024 & -0.053 & 0.054 & -0.014 & 0.014 & 5.17 & 0.035 \\ 
    & 16200 & -0.015 & 0.015 & -0.056 & 0.056 & -0.009 & 0.009 & 5.25 & 0.037 \\ 
   \hline
\end{tabular}
\caption{Simulation results in the setting \eqref{eq:cont_setting}, with conditional distribution of $W_i \cond X_i$ as in
\eqref{eq:cont_W1} (setup 1) and \eqref{eq:cont_W2} (setup 2). We report bias and root-mean squared error for the
average derivate $\theta$ based on the regression estimator \smash{$\htheta_{reg}$}, the weighted estimator
\smash{$\htheta_{weighted}$}, and the doubly robust estimator \smash{$\htheta_{DR}$}. The root mean-squared
standardized error $S$ captures the asymptotic behavior of standard Gaussian confidence intervals for $\theta$
based on \smash{$\htheta_{DR}$}. Finally, the last column reports policy value obtained by learning with doubly robust
scores over the class $\Pi$ of depth-2 trees.}
\label{tab:continuous}
\end{table}

Table \ref{tab:continuous} shows results for both average derivative estimation as described above, and for
policy learning with doubly robust scores. For policy learning, we use a cost of
treatment parameter $C = 0.2$. First, encouragingly, we see that
the doubly robust estimator of the average derivative, $\htheta_{DR}$, converges with sample size $n$,
and that the value of our learned policies improves with $n$.
Furthermore, we see that the doubly robust estimator out-performs the pure regression adjustment and weighting
estimators here. However, even the doubly-robust estimator is still bias-dominated here, and the root-mean squared standardized
error $S$ is much bigger than $1$ in all considered settings---especially the challenging ones with a non-Gaussian distribution of $W_i \cond X_i$.
This suggests that the simulation problem considered here is a difficult non-parametric problem where semiparametric
efficiency asymptotics kick in slowly at best. It is plausible that a more carefully tailored estimator of the weighting function
$d/dw \log f(x, w)$ following the lines of, e.g., \citet*{chernozhukov2018double} or \citet{hirshberg2018balancing} could
improve performance here.

\section{Proofs}
\label{app:proofs}

\subsection{Proof of Lemma \ref{lemm:rademacher}}

Our proof of this result follows the outline of the classical chaining argument of \citet{dudley1967sizes},
whereby we construct a sequence of approximating sets of increasing precision for \smash{$\tA_n(\pi)$}
with $\pi \in \Pi_n^\lambda$, and then use finite sample concentration inequalities to establish
the behavior of \smash{$\tA_n(\pi)$} over this approximation set. The improvements in our
results relative to existing bounds described in the body of the text come from
a careful construction of approximating sets targeted to the
problem of doubly robust policy evaluation---for example, our use of chaining with respect to the random
distance measure defined in \eqref{eq:D2}---and the use of sharp concentration inequalities.

Given these preliminaries, we start by defining
the conditional 2-norm distance between two policies $\pi_1, \, \pi_2$ as
\begin{equation}
\label{eq:D2}
D_n^2\p{\pi_1, \, \pi_2} =\sum_{i = 1}^n \Gamma_i^2 \p{\pi_1(X_i) - \pi_2(X_i)}^2 \, \big/ \, \sum_{i = 1}^n \Gamma_i^2, 
\end{equation}
and let $N_{D_n}(\varepsilon, \, \Pi_n^\lambda, \, \cb{X_i, \, \Gamma_i})$ be the $\varepsilon$-covering number in this distance.
To bound $N_{D_n}$, imagine creating another sample \smash{$\cb{X'_j}_{j = 1}^m$}, with $X'_j$ contained in
the support of \smash{$\cb{X_i}_{i = 1}^n$}, such that
$$\abs{\abs{\cb{j \in 1, \, ..., \, m : X'_j = X_i}} - m \, \Gamma_i^2 / \sum_{j = 1}^n \Gamma_j^2} \leq 1. $$
We immediately see that, for any two policies $\pi_1$ and $\pi_2$,
$$ \frac{1}{m} \sum_{j = 1}^m 1\p{\cb{\pi_1(X_j') \neq \pi_2(X_j')}} =  D_n^2\p{\pi_1, \, \pi_2} + \oo\p{\frac{1}{m}}. $$
Moreover, recall that the Hamming covering number $N_H$ as used in \eqref{eq:entropy_decay} does not
depend on sample size, so we can without reservations make $m$ arbitrarily large,
and conclude that
\begin{equation}
\label{eq:ND_bound}
N_{D_n}\p{\varepsilon, \, \Pi_n, \, \cb{X_i,\, \Gamma_i}} \leq N_H\p{\varepsilon^2, \, \Pi_n}.
\end{equation}
In other words, we have found that we can bound the ${D_n}$-entropy of $\Pi_n$ with
respect to its distribution-independent Hamming entropy which is controlled via Assumption \ref{assu:VC}.

Our proof strategy involves a chaining argument with respect to $D_n$. The lemma below describes the
chaining that we use in our argument; we defer the proof of Lemma \ref{lemm:chaining} to the end of this section.

\begin{lemm}
\label{lemm:chaining}
For any $J \geq 1$, there exists a chain of approximators $\Psi_j : \Pi_n^\lambda \rightarrow \Pi_n^\lambda$
for $j = 1, \, ..., \, J$, such that the following properties hold for all values of $j = 1,\, ..., \, J$
(we use the notational shorthand $\Psi_{J+1}(\pi) = \pi$):
\begin{itemize}
\item The approximation is accurate, i.e., ${D_n}(\Psi_j(\pi), \, \Psi_{j+1}(\pi)) \leq 2^{-j}$ for all $\pi \in \Pi_n^\lambda$;
\item There is no branching, such that $\Psi_j(\pi) = \Psi_j\p{\Psi_{j+1}(\pi)}$ for all $\pi \in \Pi_n^\lambda$; and
\item The set $\Pi_n^\lambda(j) := \cb{\Psi_j(\pi) : \pi \in \Pi_n^\lambda}$ of $j$-th order approximating policies
has cardinality at most $N_{D_n}(2^{-(j+1)}, \, \Pi_n, \, \cb{X_i, \, \Gamma_i})$.
\end{itemize}
\end{lemm}

We now move to our main task, i.e., bounding the Rademacher complexity $\EE{\rr_n(\Pi_n^\lambda)}$.
In order to do so, we use a two-step strategy.
We first prove the following weaker result below, with a bound that depends only on the worst-case variance
$S_n$ rather than the slice-adapted variance $S_n^\lambda \leq S_n$. We then use this bound to sharpen our
argument and prove the desired bound \eqref{eq:rad_bound}.

\begin{lemm}
\label{lemm:rademacher_weak}
Under the conditions of Lemma \ref{lemm:rademacher} and for any $\lambda$,
\begin{equation}
\label{eq:rad_bound_weak}
\limsup_{n \rightarrow \infty} \ \EE{\rr_n\p{\Pi_n^\lambda}}  \, \bigg/\, \sqrt{ \frac{S_n \VC(\Pi_n)}{n}}  \, \leq \, 52.
\end{equation} 
\end{lemm}

\begin{proof}

To start, it is helpful to decompose the random variable into several parts
using the chaining established in Lemma \ref{lemm:chaining}.
In doing so, the following thresholds play a key role:
\begin{equation}
\label{eq:Jn}
\begin{split}
J_0 := 1, \ \ \ 
J(n) := \left \lfloor \log_2(n) \p{3 - 2\beta } \, \big/\, 8 \right\rfloor, \ \eqand \
J_+(n) := \left \lfloor \log_2(n) \p{1 - \beta  } \right\rfloor.
\end{split}
\end{equation}
We then apply Lemma \ref{lemm:chaining} to create a chain with $J = J_+(n)$ terms and note that
\begin{equation}
\label{eq:pi_seq}
\begin{split}
&\frac{1}{n}\sum_{i = 1}^n \xi_i \Gamma_i  \p{2\pi(X_i) - 1}
=\frac{1}{n}\sum_{i = 1}^n \xi_i \Gamma_i  \p{2\Psi_{J_0}(\pi)(X_i) - 1} \\
&\ \ \ \ \ \
+   \sum_{j = J_0 + 1}^{J(n)} \frac{2}{n}\sum_{i = 1}^n \xi_i \Gamma_i  \p{\Psi_j(\pi)(X_i) - \Psi_{j-1}(\pi)(X_i)} \\
&\ \ \ \ \ \
+  \sum_{j = J(n) + 1}^{J_+(n)} \frac{2}{n}\sum_{i = 1}^n    \xi_i\Gamma_i \p{\Psi_j(\pi)(X_i) - \Psi_{j-1}(\pi)(X_i)} \\
&\ \ \ \ \ \
+  \frac{2}{n}\sum_{i = 1}^n \xi_i \Gamma_i  \p{\pi(X_i) - \Psi_{J_+(n)}(\pi)(X_i)},
\end{split}
\end{equation}
for any $\pi \in \Pi_n^\lambda$.
Note that, for now, the first threshold $J_0 = 1$ is trivial; however, once we want to prove the stronger bound
\eqref{eq:rad_bound} instead of \eqref{eq:rad_bound_weak} we will need a more careful choice of
$J_0$, so we already introduce this flexibility now for notational consistency. 

We now proceed to successively control the $1/\sqrt{n}$-scale behavior of all four terms above,
uniformly over all $\pi \in \Pi_n^\lambda$. The result will be that the first term
can be characterized directly via Bernstein's inequality; the second term is controlled to $1/\sqrt{n}$-scale
by chaining; the third term is shown to stochastically vanish at $1/\sqrt{n}$-scale by chaining; and
the last term is shown to deterministically vanish at $1/\sqrt{n}$-scale.

Before embarking on this task, we recall Bernstein's inequality, which will be frequently used throughout the proof:
\begin{equation}
\label{eq:bern_def}
\PP{\frac{1}{\sqrt{n}} \abs{\sum_{i = 1}^n U_i} \geq t} \leq
2\exp\sqb{\frac{-t^2}{2} \,\bigg/\, \p{\frac{1}{n} \sum_{i = 1}^n \EE{U_i^2} + \frac{Mt}{3\sqrt{n}}}},
\end{equation}
for any independent, mean-zero variables $U_i$ with $\abs{U_i} \leq M$, and any constant $t > 0$.
To make use of this inequality, it is helpful to restrict ourselves to a study of
$\rr_n(\Pi_n^\lambda)$ on the event
\begin{equation}
\label{eq:clipping}
\bb_n = \cb{M_n \leq n^{\frac{1 - 2\beta}{16}} \ \eqand \ 
\hVar{\p{2\pi(X_i) - 1}\Gamma_i} \geq \frac{s^2}{2}
\ \text{ for all } \pi \in \Pi_n^\lambda(J_0)},
\end{equation}
where $M_n = \max_{i = 1, \, ..., \, n}\cb{\abs{\Gamma_i}}$ and
$0 < \beta < 1/2$ is the constant from Assumption \ref{assu:VC}.
Recall that, by assumption, $\Gamma_i$ is sub-Gaussian and $\Var{\Gamma_i \cond X_i} > s^2$, and so
a simple calculation can be used to check that $\limn \PP{\bb_n} = 1$ and furthermore
\begin{equation}
\label{eq:clipping_2}
\limn \sqrt{n}\p{\EE{\rr_n(\Pi_n^\lambda)} - \EE{\rr_n\p{\Pi_n^\lambda} 1\p{\bb_n}}} = 0.
\end{equation}
Thus, for the rest of this proof, we will assume that the event $\bb_n$ has occurred when convenient.

\paragraph{First Term}

Because the chaining created in Lemma \ref{lemm:chaining} has no branching, we see that
\begin{equation}
\label{eq:J0_count}
\begin{split}
&\sup\cb{\frac{1}{n}\sum_{i = 1}^n \xi_i \Gamma_i  (2\Psi_{J_0}(\pi)(X_i) - 1)  : \pi \in \Pi_n^\lambda} \\
&\ \ \ \ \ \ \ \ \ \ = \sup\cb{\frac{1}{n}\sum_{i = 1}^n \xi_i \Gamma_i  (2\pi(X_i) - 1)  : \pi \in \Pi_n^\lambda(J_0)}. 
\end{split}
\end{equation}
Then, applying a union bound with Bernstein's inequality \eqref{eq:bern_def} on the event $\bb_n$ in \eqref{eq:clipping}, we
see that, for all large enough $n$ and all $t \leq  2 \hS^{0.5} \sqrt{\log(n) + \log\p{2\abs{\Pi_n^\lambda(J_0)}}}$
\begin{equation}
\label{eq:bern_first_weak}
\begin{split}
&1\p{\bb_n} \PP{\sqrt{n} \sup\cb{\frac{1}{n}\sum_{i = 1}^n \xi_i \Gamma_i  (2\pi(X_i) - 1)  : \pi \in \Pi_n^\lambda(J_0)} \geq t  \cond \cb{X_i,\, \Gamma_i}} \\
&\ \ \ \ \ \ \ \ \ \ \leq 2\abs{\Pi_n^\lambda(J_0)} \exp\sqb{-\frac{t^2}{2} \,\bigg/\, \p{\hS + t \, n^{-\frac{7 +2\beta}{16}} \,/\, 3}} \\
&\ \ \ \ \ \ \ \ \ \ \leq 2\abs{\Pi_n^\lambda(J_0)} \exp\sqb{-\frac{t^2}{4 \hS}},
\end{split}
\end{equation}
where $\hS = \sum_{i = 1}^n \Gamma_i^2/n$.
Now, to bound expectations, we note the following fact: If a non-negative random variable satisfies $X \leq c_k$ with probability $1 - 2^{-k}$
for all $k = 1,\, 2, \, \ldots$, then $\EE{X} \leq \sum_{k = 1}^\infty 2^{-k} c_k$.
Thus, applying the above bound for the choice
$$ t_k =  2 \, \hS^{0.5} \sqrt{\min\cb{k \log(2), \, \log(n)} + \log\p{2\abs{\Pi_n^\lambda(J_0)}}}, \ \ \ k = 1, \, 2, \, ..., \, \lceil \log(n) / \log(2) \rceil $$ 
we then find that (the last term corresponds to a loose $\max{\abs{\Gamma_i}}/n$ when all events fail)
\begin{equation}
\label{eq:first_count_weak}
\begin{split}
&1\p{\bb_n} \EE{\sqrt{n} \sup\cb{\frac{1}{n}\sum_{i = 1}^n \xi_i \Gamma_i  (2\pi(X_i) - 1)  : \pi \in \Pi_n^\lambda(J_0)} \cond \cb{X_i,\, \Gamma_i}} \\
&\ \ \ \ \ \ \ \ \ \ \leq 2 \, \hS^{0.5} \p{\sqrt{\log\abs{\Pi_n^\lambda(J_0)}} +  \sum_{k = 1}^\infty 2^{-k} \sqrt{(k + 1)\log(2)}} + n^{-\frac{7+2\beta}{16}}\\
&\ \ \ \ \ \ \ \ \ \ \leq 2 \, \hS^{0.5} \p{\sqrt{\log N_{H}\p{1/16, \, \Pi_n}} +  1.5} + n^{-\frac{7+2\beta}{16}}\\
&\ \ \ \ \ \ \ \ \ \ \leq 2 \, \hS^{0.5} \p{\sqrt{5 \log(16) \VC(\Pi_n)} +  1.5} + n^{-\frac{7+2\beta}{16}} \leq 11 \sqrt{\hS \, \VC(\Pi_n)} + n^{-\frac{7+2\beta}{16}},
\end{split}
\end{equation}
where for the third line we used Lemma \ref{lemm:chaining} and \eqref{eq:ND_bound} whereas for the
last line we used Assumption \ref{assu:VC} together with \eqref{eq:haussler}. Finally, noting that
\begin{equation}
\label{eq:SNconc}
\EE{\sqrt{\hS}} \leq \sqrt{S_n}
\end{equation}
by concavity of the square-root function, we see that
\begin{equation}
\label{eq:first_term_weak}
 \limsup_{n \rightarrow \infty} \ \EE{1\p{\bb_n} \sqrt{\frac{n}{S_n \VC(\Pi_n)}} \sup\cb{\frac{1}{n}\sum_{i = 1}^n \xi_i \Gamma_i  (2\pi(X_i) - 1)  : \pi \in \Pi_n^\lambda(J_0)}} \leq 11.
\end{equation}

\paragraph{Second Term}

First, we check that, for any choice of $\pi \in \Pi_n^\lambda$, $j = 1, \, ..., \, J$ and $t > 0$, we have
\begin{equation}
\label{eq:bernstein_applied}
\begin{split}
&\PP{\abs{\frac{1}{\sqrt{n}} \sum_{i = 1}^n \Gamma_i \xi_i \p{\Psi_j(\pi)(X_i) - \Psi_{j+1}(\pi)(X_i)}}
\geq t\, 2^{-j} \sqrt{\hS} \cond \cb{X_i, \, \Gamma_i}}  \\
&\ \ \ \ \ \ \ \ \ \ \ \leq 2 \exp\sqb{\frac{-t^2}{2} \p{1 + \frac{1}{3}\,\frac{M_nt\,2^j}{\sqrt{n\hS}}}^{-1}},
\end{split}
\end{equation}
where $\hS = \sum_{i = 1}^n \Gamma_i^2 / n$, $M_n = \max\cb{\abs{\Gamma_i} : 1 \leq i \leq n}$.
This can be verified using Bernstein's inequality \eqref{eq:bern_def},
which establishes that, for any choice of $t > 0$, $\pi \in \Pi_n^\lambda$ and $j = 1, \, 2, \, ..., \, J$,
\begin{align*}
&\PP{\abs{\frac{1}{\sqrt{n}} \sum_{i = 1}^n \Gamma_i \xi_i \p{\Psi_j(\pi)(X_i) - \Psi_{j+1}(\pi)(X_i)}}
\geq t\, 2^{-j} \sqrt{\hS} \cond \cb{X_i, \, \Gamma_i}} \\
&\ \  \leq  2\exp\sqb{\frac{-t^2 4^{- j}\hS}{2} \,\bigg/\, \p{\frac{1}{n} \sum_{i = 1}^n \Gamma_i^2 1\p{\cb{\Psi_j(\pi)(X_i) \neq \Psi_{j+1}(\pi)(X_i)}} + \frac{M_nt \,2^{-j} \sqrt{\hS}}{3\sqrt{n}}}} \\
&\ \  =  2\exp\sqb{\frac{-t^2 }{2} 4^{- j}\hS \,\bigg/\, \p{D_n^2\p{\Psi_j(\pi), \, \Psi_{j+1}(\pi)} \hS + \frac{M_nt \,2^{-j} \sqrt{\hS}}{3\sqrt{n}}}}.
\end{align*}
Finally recall that, by Lemma \ref{lemm:chaining}, \smash{$D_n^2\p{\Psi_j(\pi), \, \Psi_{j+1}(\pi)} \leq 4^{-j}$}; thus
$$ 4^{- j}\hS \,\bigg/\, \p{D_n^2\p{\Psi_j(\pi), \, \Psi_{j+1}(\pi)} \hS + \frac{M_nt \,2^{-j} \sqrt{\hS}}{3\sqrt{n}}} \geq \p{1 + \frac{1}{3} \frac{M_n t 2^j}{\sqrt{n \hS}}}^{-1}, $$
and so \eqref{eq:bernstein_applied} follows.

Now, or every $j \geq J_0$ and $\delta > 1/(2n)$, define the event
\begin{equation}
\label{eq:event}
\begin{split}
&\mathcal{E}_{j, \, \delta} := \cb{ \sup_{\pi \in \Pi_n^\lambda} \abs{\frac{1}{\sqrt{n}} \sum_{i = 1}^n
\Gamma_i \xi_i \p{\Psi_j(\pi)(X_i) - \Psi_{j+1}(\pi)(X_i)}} \geq  2^{-j} t_{j, \, \delta} \sqrt{\hS}}\\
&t_{j, \, \delta} := 2\sqrt{7 (j + 2) \VC(\Pi_n) + \log\p{\frac{2j^2}{\delta}}}.
\end{split}
\end{equation}
By \eqref{eq:bernstein_applied}, we immediately see that
\begin{equation}
\label{eq:step_tail}
 \PP{\mathcal{E}_{j, \, \delta}  \cond \cb{X_i, \, \Gamma_i}} \leq 2 \abs{\Pi_n^\lambda\p{j+1}} \exp\sqb{\frac{-t_{j, \, \delta}^2}{2} \p{1 + \frac{1}{3}\,\frac{M_nt_{j, \, \delta}\,2^j}{\sqrt{n\hS}}}^{-1}}.
\end{equation}
By invoking Assumption \ref{assu:VC}, Lemma \ref{lemm:chaining} and \eqref{eq:ND_bound}
along with the fact that $5\log(4) < 7$, we see that
\begin{equation}
\log\p{\abs{\Pi_n^\lambda\p{j+1}}} \leq \log\p{N_H\p{4^{-(j+2)}, \ \Pi_n}} \leq 7 (j + 2) \VC(\Pi_n).
\end{equation}
Moreover, on the event $\bb_n$ from \eqref{eq:clipping}
and recalling Assumption \ref{assu:VC} along with the definition of $J(n)$,
we see that
\begin{align*}
\frac{1}{3}\,\frac{M_nt_{j, \, \delta}\,2^j}{\sqrt{n\hS}}
&\leq \frac{2}{3} \frac{n^{\frac{1 - 2\beta}{16}} \sqrt{7(J(n) + 2) \VC(\Pi_n) + \log(2 n J(n)^2)} 2^{J(n)}}{\sqrt{n s^2 / 2}} \\
&= \exp\sqb{\log(n)\p{\frac{1 - 2\beta}{16} + \frac{\beta}{2} + \frac{3 - 2\beta}{8} - \frac{1}{2}}} \cdot \text{polylog}(n) \\
&= n^{\frac{2\beta - 1}{16}} \cdot \text{polylog}(n) \leq 1
\end{align*}
for large enough values of $n$, 
simultaneously for all $j \leq J(n)$ and $\delta \geq 1/(2n)$, because $\beta < 1/2$.
Thus, for large enough values of $n$, the bound \eqref{eq:step_tail} simplifies dramatically, and we get
\begin{equation}
 1\p{\bb_n} \PP{\mathcal{E}_{j, \, n} \cond \cb{X_i, \, \Gamma_i}} \leq \frac{\delta}{j^2}.
\end{equation}
Applying this bound simultaneously to $j = J_0, \, ..., \, J(n) - 1$:
\begin{equation}
\begin{split}
  1\p{\bb_n} \PP{\bigcup_{j = J_0}^{J(n) - 1} \mathcal{E}_{j, \, n} \cond \cb{X_i, \, \Gamma_i}}
 &\leq \sum_{j = J_0}^{J(n) - 1} \frac{\delta}{j^2} \leq 2\delta.
 \end{split}
\end{equation}
Thus, for large enough $n$, we can directly verify that, with probability at least $1 - 2\delta$,
\begin{equation*}
\begin{split}
&\sqrt{n}1\p{\bb_n} \sup_{\pi \in \Pi_n^\lambda} \abs{\frac{2}{n} \sum_{i=1}^n \Gamma_i \xi_i \sum_{j = J_0}^{J(n) - 1} \p{\Psi_{j+1}(\pi) - \Psi_{j}(\pi)}(X_i)} \\
&\ \ \ \ \ \ \ \leq 4\sqrt{\hS} \sum_{j = J_0}^{J(n) - 1} 2^{-j} \sqrt{7 (j + 2) \VC(\Pi_n) + \log\p{\frac{2j^2}{\delta}}} \\
&\ \ \ \ \ \ \ \leq 4 \sqrt{\hS} \p{\sqrt{7\VC(\Pi_n)} \sum_{j = J_0}^{J(n) - 1} 2^{-j} \sqrt{j+2} + \sum_{j = J_0}^{J(n) - 1} 2^{-j}  \sqrt{\log\p{2j^2}} + 2^{1-J_0}\sqrt{\log\p{\delta^{-1}}}}.
\end{split}
\end{equation*}
Moreover, we can check by calculus that, for all $J_0 \geq 2$,
\begin{align*}
&\sum_{j = J_0}^{J(n) - 1} 2^{-j} \sqrt{j+2}  \leq 2^{- J_0}  \sum_{j = 0}^\infty 2^{-j} \p{\sqrt{J_0} + \frac{j + 2}{2 \sqrt{J_0}}} 
= 2 \times 2^{- J_0} \sqrt{J_0} + 3 \times 2^{-J_0}, \\
& \sum_{j = J_0}^{J(n) - 1} 2^{-j}  \sqrt{\log\p{2j^2}} \leq 2^{ - J_0}  \sum_{j = 0}^\infty 2^{-j} \p{\sqrt{\log\p{2J_0^2}} + \frac{2\log(J_0 + j) - 2\log(J_0)}{2 \sqrt{\log\p{2J_0^2}}}} \\
&\ \ \ \ \ \ \ \ \ \ \ \ \ \ \ \ \ \leq 2^{ - J_0}  \sum_{j = 0}^\infty 2^{-j} \p{\sqrt{\log\p{2J_0^2}} + \frac{j}{J_0 \sqrt{\log\p{2J_0^2}}}} \\
&\ \ \ \ \ \ \ \ \ \ \ \ \ \ \ \ \ = 2 \times 2^{- J_0} \p{\sqrt{\log\p{2J_0^2}} + \frac{1}{J_0\sqrt{\log\p{2J_0^2}}}} \leq 4 \times 2^{- J_0}  \sqrt{J_0};
\end{align*}
moreover, the same final upper bounds can be verified directly for $J_0 = 1$.
Thus the above expression can further be bounded by
$$ \ldots \leq 4\sqrt{\hS} 2^{-J_0} \p{\sqrt{7 \VC(\Pi_n)} \p{2 \sqrt{J_0} + 3} + 4 \sqrt{J_0}  + 2 \sqrt{\log(\delta^{-1})}}. $$
Next, we bound expectations as in \eqref{eq:first_count_weak}, and apply the above bound
separately for the sequences $2\delta = \max\cb{2^{-k}, \, 1/n}$ for $k = 1, \, 2, \, ...$
to show that, again for large enough $n$, 
\begin{equation}
\label{eq:chain_tail}
\begin{split}
&\sqrt{n} \EE{1\p{\bb_n} \sup_{\pi \in \Pi_n^\lambda} \abs{\frac{2}{n} \sum_{i=1}^n \Gamma_i \xi_i \sum_{j = J_0}^{J(n) - 1} \p{\Psi_{j+1}(\pi) - \Psi_{j}(\pi)}(X_i)} } \\
&\ \ \ \ \ \ \ \ \leq 4 \times 2^{-J_0} \p{\sqrt{7 \VC(\Pi_n)} \p{2 \sqrt{J_0} + 3} + 4 \sqrt{J_0}  + 2 \sum_{k = 1}^\infty 2^{-k}\sqrt{(k+1)\log(2)}} \EE{\sqrt{\hS}} \\
&\ \ \ \ \ \ \ \ \leq 2^{-J_0}\sqrt{\VC(\Pi_n)}  \p{38 \sqrt{J_0} + 44} \EE{\sqrt{\hS}} \leq 2^{-J_0} \sqrt{S_n \VC(\Pi_n)}  \p{38 \sqrt{J_0} + 44},
\end{split}
\end{equation}
where we note that the contribution of terms on the residual with-probability-$1/n$ scale
as $M_n / n \ll 1/\sqrt{n}$ on $\bb_n$ \eqref{eq:clipping}, and
for the last inequality we also use \eqref{eq:SNconc}.
We thus conclude that
\begin{equation}
\label{eq:second_term_final_weak}
\begin{split}
&\limsup_{n \rightarrow \infty} \sqrt{\frac{n}{S_n \VC(\Pi_n)}} \EE{1\p{\bb_n} \sup_{\pi \in \Pi_n^\lambda} \abs{\frac{1}{n} \sum_{i=1}^n \Gamma_i \xi_i \sum_{j = J_0}^{J(n) - 1} \p{\Psi_{j+1}(\pi) - \Psi_{j}(\pi)}(X_i)} } \leq 41,
\end{split}
\end{equation}
recalling our choice of $J_0 = 1$ from \eqref{eq:Jn}.

\paragraph{Third Term}

We now verify that terms $\Psi_j(\pi)(X_i) - \Psi_{j+1}(\pi)(X_i)$ in \eqref{eq:pi_seq}
with $J(n) \leq j < J_+(n)$ are asymptotically negligible.
To do so, we collapse all approximating policies with $J(n) \leq j < J_+(n)$, and directly
compare $\Psi_{J(n)}(\pi)$ to $\Psi_{J_+(n)}(\pi)$. Because of our ``no branching'' construction,
we know that $\Psi_{J(n)}(\pi) = \Psi_{J(n)}(\Psi_{J_+(n)}(\pi))$ for all policies $\pi \in \Pi_n^\lambda$, and so
\begin{align*}
&\PP{\sup\cb{\abs{\frac{1}{\sqrt{n}} \sum_{i = 1}^n \Gamma_i \xi_i \p{\Psi_{J(n)}(\pi)(X_i) - \Psi_{J_+(n)}(\pi)(X_i)}} : \pi \in \Pi_n^\lambda}
\geq  2 \times  t\, 2^{-J(n)} \sqrt{\hS}}  \\
&\ \ \ \ \ \ =\PP{\sup\cb{\abs{\frac{1}{\sqrt{n}} \sum_{i = 1}^n \Gamma_i \xi_i \p{\Psi_{J(n)}(\pi)(X_i) - \pi(X_i)}} : \pi \in \Pi_n^\lambda\p{J_+(n)}}
\geq 2 \times  t\, 2^{-J(n)} \sqrt{\hS}} \\ 
&\ \ \ \ \ \ \leq 2 \abs{\Pi_n^\lambda\p{J_+(n)}}\exp\sqb{\frac{-t^2}{2}\p{1 + \frac{1}{6}\,\frac{M_nt\,2^{J(n)}}{\sqrt{n\hS}}}^{-1}},
\end{align*}
where the last inequality follows from Bernstein's inequality using exactly the same arguments as those used to establish \eqref{eq:bernstein_applied}.
By Lemma \ref{lemm:chaining}, Assumption \ref{assu:VC} and \eqref{eq:haussler}, we get
\begin{equation}
\label{eq:card_bound}
\begin{split}
&\log \abs{\Pi_n^\lambda\p{J_+(n)}} \leq \log  N_{D_n}\p{2^{-(J_+(n) + 1)}, \, \Pi_n, \, \cb{X_i, \, \Gamma_i}} \\
&\ \ \ \ \ \ \ \  \leq \log  N_H\p{4^{-(J_+(n) + 1)}, \, \Pi_n} \leq  5\log(4)  (J_+(n) + 1) n^\beta. 
\end{split}
\end{equation}
The next step is to plug \smash{$t^2 = 4^{J(n)}  n^{(2\beta - 1)/4} /\hS$} into the previous bound.
Given this choice along with Assumption \ref{assu:VC} and \eqref{eq:Jn} we see that, on event $\bb_n$ from \eqref{eq:clipping},
$$ {t 2^{J(n)}}\,/{\sqrt{n}} \geq  \times 2^{2J(n)} n^{\frac{2\beta - 5}{8}} n^{\frac{-1 + 2\beta}{16}} \geq n^{\frac{1 - 2\beta}{16}} / 4 $$
which grows with $n$, and so the bound simplifies on event $\bb_n$ and for large enough $n$:
\begin{align*}
&\PP{1\p{\bb_n} \Delta_{mid}\p{\Pi_n^\lambda} \geq 2 n^{\frac{2\beta - 1}{8}}} 
\leq 1\p{\bb_n}  2 \abs{\Pi_n^\lambda\p{J_+(n)}}\exp\sqb{\frac{-(3/2)t\sqrt{n\hS}}{2^{J(n)} \max\cb{M_n, \, 1}}}  \\
&\ \ \ \ \ \ \leq  2 \exp\sqb{\sqrt{n} \p{5\log(4)  (J_+(n) + 1) n^{\beta - 1/2} - \frac{3}{2} n^{\frac{6\beta - 3}{16}}} }, \ \where \\
&\Delta_{mid}\p{\Pi_n^\lambda} = \sup\cb{\abs{\frac{1}{\sqrt{n}} \sum_{i = 1}^n \Gamma_i \xi_i \p{\Psi_{J(n)}(\pi)(X_i) - \Psi_{J_+(n)}(X_i)}} : \pi \in \Pi_n^\lambda}.
\end{align*}
Thus, noting that $\beta < 1/2$,  we see that
\begin{align*}
&\limsup_{n \rightarrow \infty} n^{\frac{5 + 6\beta}{16}} \log\p{ \PP{1\p{\bb_n} \Delta_{mid}\p{\Pi_n^\lambda}}
\geq 2 n^{\frac{2\beta - 1}{8}}}  \leq -\frac{3}{2}.
\end{align*}
Meanwhile, we also know that $1\p{\bb_n} \Delta_{mid}\p{\Pi_n^\lambda} / \sqrt{n} \leq n^{(1 - 2\beta)/16}$, and so we conclude that
$$ \limn \EE{\sup\cb{\abs{\frac{1}{\sqrt{n}} \sum_{i = 1}^n \Gamma_i \xi_i \p{\Psi_{J(n)}(\pi)(X_i) - \Psi_{J_+(n)}(X_i)}} : \pi \in \Pi_n^\lambda}} = 0, $$
meaning that the third group of terms in the chaining \eqref{eq:pi_seq} in fact do not contribute to the first-order behavior of the Rademacher complexity.

\paragraph{Fourth Term}

Finally, the last term in \eqref{eq:pi_seq} can be shown to vanish at $1/\sqrt{n}$-scale determinisitically. By Cauchy-Schwarz,
\begin{align*}
&\abs{\frac{1}{n} \sum_{i = 1}^n \Gamma_i \xi_i \p{\pi(X_i) - \Psi_{J_+(n)}(\pi)(X_i)}}
\leq \sqrt{\frac{1}{n} \sum_{i = 1}^n \Gamma_i^2 \p{\pi(X_i) - \Psi_{J_+(n)}(\pi)(X_i)}^2} \\
&\ \ \ \ \ \ \ \ = {D_n}\p{\pi, \, \Psi_{J_+(n)}(\pi)} \sqrt{\hS}
\leq 2^{-J_+(n)} \sqrt{\hS}.
\end{align*}
Furthermore, recalling the definition of $J_+(n)$ from \eqref{eq:Jn} and
on the event where $M_n$ is controlled as in \eqref{eq:clipping},
$$\limn  \sqrt{n} 2^{-J_+(n)} \sqrt{\hS} \leq 2  \sqrt{n} n^{\beta - 1} n^{\frac{1 - 2\beta}{16}}  = n^{\frac{14\beta - 7}{16}} = 0,$$
because  $\beta < 1/2$ by Assumption \ref{assu:VC}.

\paragraph{Wrapping Up Lemma \ref{lemm:rademacher_weak}}

Combining \eqref{eq:first_term_weak} with \eqref{eq:second_term_final_weak} with our above results showing
that the third and fourth terms in \eqref{eq:pi_seq} are asymptotically negligible, we recover \eqref{eq:rad_bound_weak}.
\end{proof}


We now turn to proving Lemma \ref{lemm:rademacher} itself, and specifically the bound \eqref{eq:rad_bound}.
In doing so, we follow the proof of Lemma \ref{lemm:rademacher_weak} closely, but with slightly stronger concentration
bounds that are unlocked by the result we already have in Lemma \ref{lemm:rademacher_weak}. We also replace the
choice $J_0 = 1$ in \eqref{eq:Jn} with
\begin{equation}
\label{eq:J0_new}
J_0 := 9 + \left\lfloor \log_4 \p{ {S_n}\,\big/\,{S_n^\lambda}} \right\rfloor.
\end{equation}
In the resulting new decomposition \eqref{eq:pi_seq}, we note that the third and fourth terms are still vanishing
at the $1/\sqrt{n}$-scale, so we do not need to revisit those. Thus, our only task is to sharpen our bounds on the
first and second terms.

The main additional work we need to do is in bounding the first term. Starting from \eqref{eq:J0_count}
we note that, because the $\xi_i$ are all mean-zero,
\begin{equation}
\label{eq:centering}
\begin{split}
& \EE{\sup\cb{\frac{1}{n}\sum_{i = 1}^n \xi_i \Gamma_i  (2\pi(X_i) - 1) : \pi \in \Pi_n^\lambda(J_0)}} \\
&\ \ \ \ \ \ \ \ \ \ = \EE{\sup\cb{\frac{1}{n}\sum_{i = 1}^n \xi_i \p{\Gamma_i  (2\pi(X_i) - 1)   - A_n^*} : \pi \in \Pi_n^\lambda(J_0)}},
\end{split}
\end{equation}
where $A_n^* = \sup\cb{A_n(\pi) : \pi \in \Pi_n^\lambda}$. Then, applying Bernstein's inequality as in
\eqref{eq:bern_first_weak}, we get that
for all large enough $n$ and all $t \leq  2 \hS_{\max}^{0.5} \sqrt{\log(n) + \log\p{2\abs{\Pi_n^\lambda(J_0)}}}$,
\begin{equation}
\label{eq:bern_first_new}
\begin{split}
&1\p{\bb_n} \PP{\sqrt{n} \sup\cb{\frac{1}{n}\sum_{i = 1}^n \xi_i \Gamma_i  (2\pi(X_i) - 1)  : \pi \in \Pi_n^\lambda(J_0)} \geq t  \cond \cb{X_i,\, \Gamma_i}} \\
&\ \ \ \ \ \ \ \ \ \ \ \ \ \ \leq 2 \abs{\Pi_n^\lambda(J_0)} \exp\sqb{-\frac{t^2}{4 \hS_{\max}}}, \\
&\hS_{\max} := \sup\cb{ \frac{1}{n}\sum_{i = 1}^n \p{\Gamma_i  (2\pi(X_i) - 1)  - A_n^*}^2 : \pi \in \Pi_n^\lambda(J_0)}.
\end{split}
\end{equation}
Then, following \eqref{eq:first_count_weak}, we get that
\begin{equation}
\label{eq:first_count_new}
\begin{split}
&1\p{\bb_n} \EE{\sqrt{n} \sup\cb{\frac{1}{n}\sum_{i = 1}^n \xi_i \Gamma_i  (2\pi(X_i) - 1)  : \pi \in \Pi_n^\lambda(J_0)} \cond \cb{X_i,\, \Gamma_i}} \\
&\ \ \ \ \ \ \ \ \ \ \leq 2 \, \hS_{\max}^{0.5} \p{\sqrt{\log N_{H}\p{4^{-(J_0 + 1)}, \, \Pi_n}} +  1.5} + n^{-\frac{7 + 2\beta}{16}} \\
&\ \ \ \ \ \ \ \ \ \ \leq 2 \, \hS_{\max}^{0.5} \p{\sqrt{5 \log(4) \VC(\Pi_n) (J_0 + 1)} +  1.5} + n^{-\frac{7 + 2\beta}{16}} \\
&\ \ \ \ \ \ \ \ \ \  \leq 6 \sqrt{\hS_{\max} \, \VC(\Pi_n) \, \p{10 + \left\lfloor \log_4 \p{ {S_n}\,\big/\,{S_n^\lambda}} \right\rfloor}} + n^{-\frac{7 + 2\beta}{16}}.
\end{split}
\end{equation}
Now, combining the bound we already have from Lemma \ref{lemm:rademacher_weak} with the proof of Lemma \ref{lemm:coupling},
we see that under the conditions of Lemma \ref{lemm:rademacher} and provided that $S_n \VC(\Pi_n) / n \rightarrow 0$, we
have that
$$ \limsup_n \EE{\sqrt{\hS_{\max}}} \, \big/ \, \sqrt{S_n^\lambda + 4 \lambda^2}  \leq 1; $$
to check this, we also used the fact that, by \eqref{eq:SL},
\begin{equation}
\begin{split}
&\sup\cb{\EE{\p{2(\pi(X_i) - 1)\Gamma_i - A_n^*}^2} : \pi \in \Pi_n^\lambda} \\
&\ \ \ \ \ \ \ \ \ = \sup\cb{\Var{2(\pi(X_i) - 1)\Gamma_i} + \p{A_n(\pi) - A_n^*}^2 : \pi \in \Pi_n^\lambda}
\leq S_n^\lambda + 4\lambda^2.
\end{split}
\end{equation}
Thus, we conclude that
\begin{equation}
\label{eq:first_term_new}
\begin{split}
& \limsup_{n \rightarrow \infty} \EE{1\p{\bb_n} \sqrt{\frac{n}{\p{S_n^\lambda + 4\lambda^2} \VC(\Pi_n) } } \sup\cb{\frac{1}{n}\sum_{i = 1}^n \xi_i \Gamma_i  (2\pi(X_i) - 1)  : \pi \in \Pi_n^\lambda(J_0)}} \\
&\ \ \ \ \ \ \ \ \ \  \Big / \p{1 + 18\sqrt{1 +  \left\lfloor  \log_4 \p{ {S_n}\,\big/\,{S_n^\lambda}} \right\rfloor \Big/\, 9} } \leq 1.
\end{split}
\end{equation}
Meanwhile, for the second term, we proceed exactly as before up to \eqref{eq:chain_tail}. Here, however,
we invoke the new (larger) choice of $J_0$ and, noting that
$$ 2^{-8} \p{44 + 38 \sqrt{ 9 + \left\lfloor \log_4 \p{ \frac{S_n}{S_n^\lambda}} \right\rfloor}} \leq  \sqrt{ 1 + \left\lfloor \log_4 \p{ \frac{S_n}{S_n^\lambda}} \right\rfloor \Big/\, 9}, $$
we get
\begin{equation}
\label{eq:second_term_final}
\begin{split}
&\limsup_{n \rightarrow \infty} \sqrt{\frac{n}{S_n^\lambda \VC(\Pi_n)}} \EE{1\p{\bb_n} \sup_{\pi \in \Pi_n^\lambda} \abs{\frac{1}{n} \sum_{i=1}^n \Gamma_i \xi_i \sum_{j = J_0}^{J(n) - 1} \p{\Psi_{j+1}(\pi) - \Psi_{j}(\pi)}(X_i)} } \\
&\ \ \ \ \ \ \ \ \bigg / \sqrt{ 1 + \left\lfloor \log_4 \p{ \frac{S_n}{S_n^\lambda}} \right\rfloor \Big/\, 9} \leq 1.
\end{split}
\end{equation}
Finally, we establish \eqref{eq:rad_bound} by combining this bound with \eqref{eq:first_term_new},
and the fact that clipping as in \eqref{eq:clipping} has an asymptotically negligible effect.

\subsubsection*{Proof of Lemma \ref{lemm:chaining}}
We construct the chaining by backwards recursion, as follows. First, for the largest index
$J$ under consideration, we do the following:
\begin{enumerate}
\item Let $\Psi_J': \Pi_n \rightarrow \cb{\xx \rightarrow \cb{0, \, 1}}$ be an optimal $2^{-(J+1)}$ covering of $\Pi_n$, such that the
cardinality of the set $ \cb{\Psi_J'(\pi) : \pi \in \Pi_n}$ is at most $N_{D_n}\p{2^{-(J+1)}, \, \Pi_n, \, \cb{X_i,\, \Gamma_i}}$.
\item For every approximating policy $\pi' \in \cb{\Psi_J'(\pi) : \pi \in \Pi_n}$, construct a function $\text{neighbor}(\cdot)$ such
that $\text{neighbor}(\pi') \in \cb{\pi \in \Pi_n^\lambda : D_n(\pi, \, \pi') \leq 2^{-(J+1)}}$ if this set is non-empty, and
$\text{neighbor}(\pi') = \emptyset$ else.
\item Define $\Psi_J : \Pi_n^\lambda \rightarrow \Pi_n^\lambda$ via $\Psi_J(\pi) = \text{neighbor}(\Psi_j'(\pi))$.
\end{enumerate}
We can see by construction that $\Psi_J(\pi) \in \Pi_n^\lambda$ for all $\pi \in \Pi_n^\lambda$ (because no element in
$\Pi_n^\lambda$ can be mapped by $\Psi_J'$ to an element $\pi'$ with $\text{neighbor}(\pi') = \emptyset$), and that
the cardinality of the set $\Pi_n^\lambda(J) = \cb{\Psi_J(\pi) : \pi \in \Pi_n^\lambda}$ is at most $N_{D_n}\p{2^{-(J+1)}, \, \Pi_n, \, \cb{X_i,\, \Gamma_i}}$.
Furthermore, by the triangle inequality, $D_n(\Psi_J(\pi), \, \pi) \leq 2^{-J}$ for all $\pi \in \Pi_n^\lambda$.

Next, for every $1 \leq j < J$, we first define the mapping $\Psi_j$ as a $2^{-j}$-approximation of $\Pi_n^\lambda(j+1)$
using exactly the same construction as above. Thus, $\Psi_{j} : \Pi_n^\lambda(j+1) \rightarrow \Pi_n^\lambda(j+1)$,
$\Pi_n^\lambda(j) = \cb{\Psi_j(\pi) : \pi \in \Pi_n^\lambda(j+1)}$ has cardinality at most $N_{D_n}\p{2^{-(j+1)}, \, \Pi_n, \, \cb{X_i,\, \Gamma_i}}$,
and $D_n(\Psi_j(\pi), \, \pi) \leq 2^{-j}$ for all $\pi \in \Pi_n^\lambda(j+1)$.
Finally, we extend the mappings $\Psi_j$ to the whole domain $\Pi_n^\lambda$ via the relationship
$\Psi_j(\pi) = \Psi_j\p{\Psi_{j+1}(\pi)}$ for all $\pi \in \Pi_n^\lambda$. Note that this extension does not grow
the size of the set $\Pi_n^\lambda(j)$, and that the mapping $\Psi_j$ has no branching by construction.

\subsection{Proof of Corollary \ref{coro:additive}}

First, as argued by \citet{bartlett2002rademacher} in the proof of their Theorem 8,
\begin{equation}
\label{eq:rademacher_apply}
\EE{\sup\cb{\abs{\tA_n\p{\pi} - A_n\p{\pi}} : \pi \in \Pi_n^\lambda}} \leq 2\EE{\rr_n\p{\Pi_n^\lambda}},
\end{equation}
Then, to check concentration, we need to bound
\smash{$\sup_{\pi \in \Pi_n}|\tA_n\p{\pi} - A_n\p{\pi}|$} in terms of its expectation.
Recall that \smash{$\tA_n(\pi) = n^{-1} \sum \Gamma_i (2\pi(X_i) - 1)$}, and that the \smash{$\Gamma_i$}
are uniformly sub-Gaussian. Because the \smash{$\Gamma_i$} are not bounded, it is convenient to define
truncated statistics
$$ \tA_n^{(-)}(\pi) = \frac{1}{n} \sum_{i = 1}^n \Gamma_i^{(-)} (2\pi(X_i) - 1), \ \ \Gamma_i^{(-)} = \Gamma_i \, \ind\p{\cb{\abs{\Gamma_i} \leq \log(n)}}. $$
Here, we of course have that \smash{$|\Gamma_i^{(-)}| \leq \log(n)$}, and so we can
apply Talagrand's inequality as described in \citet{bousquet2002bennett} to these
truncated statistics. We see that, for any $\delta >0$, with probability at least $1 - \delta$,
\begin{align*}
&\sup_{\pi \in \Pi_n^\lambda} \abs{\tA_n^{(-)}\p{\pi} - A_n^{(-)}\p{\pi}} \leq \EE{\sup_{\pi \in \Pi_n^\lambda} \abs{\tA_n^{(-)}\p{\pi} - A_n^{(-)}\p{\pi}}} + \frac{\log(n)\log(\delta)}{3n}\\
&\ \ \ \ \ \ \ \  + \sqrt{2 \log\p{\delta^{-1}} \p{\sup_{\pi \in \Pi_n^\lambda}\Var{\tA_n(\pi)} + \frac{2\log(n)}{n}\EE{\sup_{\pi \in \Pi_n^\lambda} \abs{\tA_n^{(-)}\p{\pi} - A_n^{(-)}\p{\pi}}}}},
\end{align*}
where we used the short-hand \smash{$A_n^{(-)}\p{\pi} = \mathbb{E}[\tA_n^{(-)}\p{\pi}]$}.
Moreover, because the \smash{$\Gamma_i$} are uniformly sub-Gaussian, we can immediately verify that
$$ \EE{\abs{\sup_{\pi \in \Pi_n^\lambda} \abs{\tA_n^{(-)}\p{\pi} - A_n^{(-)}\p{\pi}} - \sup_{\pi \in \Pi_n} \abs{\tA_n\p{\pi} - A_n\p{\pi}}}} $$
decays exponentially fast in $n$; similarly,
\smash{$n\sup_{\pi \in \Pi_n^\lambda}\Var{\tA_n(\pi)} - S_n^\lambda$} also decays exponentially fast.
Using \eqref{eq:rademacher_apply} and noting that,
by Lemma \ref{lemm:rademacher} and Assumption \ref{assu:VC}, $\EE{\rr_n\p{\Pi_n^\lambda}}$ decays 
polynomially in $n$, we conclude
that with probability at least $1 - \delta$,
\begin{equation}
\label{eq:conc_1}
\begin{split}
&\sup\cb{ \abs{\tA_n\p{\pi} - A_n\p{\pi}} : \pi \in \Pi_n^\lambda} \\
&\ \ \ \ \ \ \ \ \leq \p{1 + o(1)} \p{\EE{\sup\cb{\abs{\tA_n\p{\pi} - A_n\p{\pi}} : \pi \in \Pi_n^\lambda}} + \sqrt{\frac{2 S_n^\lambda\log\p{\delta^{-1}}}{n}}},
\end{split}
\end{equation}
thus establishing our second claim.

\subsection{Proof of Lemma \ref{lemm:coupling}}

In the argument below, we omit all $n$-subscripts for readability, e.g., we write
\smash{$\hA(\pi)$} instead of  \smash{$\hA_n(\pi)$}.
For any fixed policy $\pi$, we begin by expanding out the difference of interest as
\begin{align*}
&\hA(\pi) - \tA(\pi) \\
&\ \ = \frac{1}{n}\sum_{i = 1}^n (2\pi(X_i) - 1) \p{Y_i - m(X_i, \, W_i)} \p{\hg^{(-k(i))}(X_i, \, Z_i) - g(X_i, \, Z_i)} \\
&\ \ \ \ + \frac{1}{n} \sum_{i = 1}^n (2\pi(X_i) - 1)\p{\tau_{\hatm^{(-k)}}(X_i \, W_i) - \tau_m(X_i, \, W_i) - g(X_i, \, Z_i) \p{\hatm^{(-k(i))}(X_i, \, W_i) - m(X_i, \, W_i)}} \\
&\ \ \ \ - \frac{1}{n} \sum_{i = 1}^n (2\pi(X_i) - 1)\p{\hatm^{(-k(i))}(X_i, \, W_i) - m(X_i, \, W_i)}\p{\hg^{(-k(i))}(X_i, \, Z_i) - g(X_i, \, Z_i)} .
\end{align*}
Denote these three summands by $D_1(\pi)$, $D_2(\pi)$ and $D_3(\pi)$.
We will bound all 3 summands separately.

To bound the first term, it is helpful separate out the contributions of the $K$ different folds:
\begin{equation}
\label{eq:DK}
\begin{split}
D^{(k)}_{1}(\pi) = \frac{1}{n} \sum_{\cb{i : k(i) = k}} &(2\pi(X_i) - 1) \p{Y_i - m(X_i, \, W_i)} \\
&\ \ \ \ \  \p{\hg^{(-k(i))}(X_i, \, Z_i) - g(X_i, \, Z_i)}.
\end{split}
\end{equation}
Now, because \smash{$\hg^{(-k)}(\cdot)$} was only computed using data from the
$K - 1$ folds, we can condition on the value of this function estimate to make the individual terms
in the above sum independent. Moreover, by exogeneity of the instrument and the exclusion restriction,
we see that \smash{$\EE{Y_i - m(X_i, \, W_i) \cond X_i, \, Z_i, \, \hg^{(-k(i))}(\cdot)} = 0$}, and
so the expected second moment of \smash{$D^{(k)}_{1}(\pi)$} reduces to the sum of the variances
of its constituent terms.

Next, by Assumption \ref{assu:LR}, we know that
$$ \sup_{x \in \xx} \abs{\p{\hg^{(-k)}(x, \, z) - g(x, \, z)}}\leq 1 $$
with probability tending to 1, and so the individual summands in \eqref{eq:DK}
are all $\nu$-sub-Gaussian with probability tending to 1. Then, writing
$$ V_n(k) = \EE{ \p{\hg^{(-k)}(X_i, \, Z_i) - g(X_i, \ Z_i)}^2 \Var{Y_i - m(X_i, \, W_i) \cond X_i, \, Z_i} \cond \hg^{(-k)}(\cdot)} $$
for the variance of \smash{$D^{(k)}_1(\pi)$} conditionally on the model
\smash{$\hg^{(-k)}(\cdot)$} fit on the other $K - 1$ folds,
we can apply Corollary \ref{coro:additive} to establish that
\begin{equation}
\label{eq:AK_bound}
\frac{n}{n_k} \ \EE{ \sup_{\pi \in \Pi} \abs{D^{(k)}_{1}(\pi)} \ \Bigg | \  \hg^{(-k)}(\cdot)}
 = \oo\p{ \sqrt{\VC\p{\Pi_n} \,\frac{V_n(k) }{n_k}}},
\end{equation}
where $n_k = \abs{\cb{i : k(i) = k}}$ denotes the number of observations in the $k$-th fold.
Since we compute our doubly robust scores using a finite number of evenly-sized
folds, $n_k/n \rightarrow 1/K$, we can use our risk bounds in Assumption \ref{assu:LR} 
to check that
\begin{equation}
\begin{split}
\EE{V_n(k)}
&\leq \EE{\nu^2 \, \EE{ \p{\hg^{(-k)}(X_i, \, Z_i) - g(X_i, \, Z_i)}^2 \cond \hg^{(-k)}(\cdot)}}  \\
&= \oo\p{a\p{\frac{K - 1}{K} \, n} \,n^{-\zeta_g}}.
\end{split}
\end{equation}
Then, applying \eqref{eq:AK_bound} separately to all $K$ folds and using Jensen's inequality, we find that
\begin{equation}
\EE{\sup_{\pi \in \Pi} \abs{D_1(\pi)}} = \oo\p{ \nu\sqrt{\VC\p{\Pi_n} \,\frac{a((1 - K^{-1})n)}{n^{1 + \zeta_g}}}},
\end{equation}
thus bounding the first term.

Meanwhile, recall that by the properties of our weighting function \eqref{eq:representer}, we know that
$\EE{\tau_{\tilde{m}}(X_i, \, W_i) - g(X_i, \, Z_i) \tilde{m}(X_i, \, W_i) \cond X_i} = 0$ for any conditional response function $\tilde{m}(\cdot)$,
which in particular means that, by cross-fitting,
\begin{align*}
&\mathbb{E}\Big[\tau_{\hatm^{(-k)}}(X_i, \, W_i) - \tau_m(X_i, \, W_i) \\
&\ \ \ \ \ \ \ \ \ - g(X_i, \, Z_i) \p{\hatm^{(-k(i))}(X_i, \, W_i) - m(X_i, \, W_i)} \cond X_i, \, \hatm^{-k(i)}(\cdot)\Big] = 0.
\end{align*}
Thus, by a similar argument as before, we find that
\begin{equation}
\EE{\sup_{\pi \in \Pi} \abs{D_{2}(\pi)}} = \oo\p{ \frac{1}{\eta}\sqrt{\VC\p{\Pi_n} \,\frac{a((1 - K^{-1})n)}{n^{1 + \zeta_m}}}},
\end{equation}
where $\eta$ is the uniform ``overlap'' bound on the weighting function $g(\cdot)$.

It now remains to bound the final term, $D_{3}(\pi)$. Here, we can use
the Cauchy-Schwarz inequality to verify that
\begin{align*}
&\abs{D_{3}(\pi)} = \bigg| \frac{1}{n} \sum_{i = 1}^n (2\pi(X_i) - 1) \p{\hatm^{(-k(i))}(X_i, \, W_i) - m(X_i, \, W_i)} \\
&\ \ \ \ \ \ \ \ \ \ \ \ \ \ \ \ \ \ \ \ \ \ \ \ \p{\hg^{(-k(i))}(X_i, \, Z_i) - g(X_i, \, Z_i)}\bigg|   \\
&\ \ \ \ \ \leq  \sqrt{\frac{1}{n}\sum_{i = 1}^n \p{\hatm^{(-k(i))}(X_i, \, W_i) - m(X_i, \, W_i)}^2} \\
&\ \ \ \ \ \ \ \ \ \ \ \ \ \ \ \ \ \ \ \ \ \ \ \ \sqrt{\frac{1}{n}\sum_{i = 1}^n \p{\hg^{(-k(i))}(X_i, \, Z_i) - g(X_i, \, Z_i)}^2}.
\end{align*}
This bound is deterministic and does not depend on $\pi$; thus, it also holds as a bound for the
supremum of $|D_3(\pi)|$ over all $\pi$.
Then, applying Cauchy-Schwarz again to the above product, we see that
\begin{align*}
\EE{\frac{n \, \sup_{\pi \in \Pi} \abs{D_3(\pi)}}{\abs{\cb{i : W_i = 1}}}}
&\leq \sqrt{\EE{\p{\hatm^{(-k(i))}(X_i, \, W_i) - m(X_i, \, W_i)}^2}} \\
&\ \ \ \ \ \ \ \ \ \sqrt{\EE{\p{\hg^{(-k(i))}(X_i, \, Z_i) - g(X_i, \, Z_i)}^2}} \\
&\leq {a\p{\left\lfloor\frac{K-1}{K} \,n\right \rfloor}} \,\Big/\, {\sqrt{\left\lfloor\frac{K-1}{K} \,n\right \rfloor}},
\end{align*}
The desired conclusion now follows from combining these three bounds.

\subsection{Proof of Theorem \ref{theo:lb}}

\newcommand{\acal}{\mathcal{A}}

Writing $\VC(\Pi) = d$, we know that there exists a collection of $d$ non-overlapping sets $\acal_j$
for $j = 1, \, ..., \, d$ such that $\Pi$ shatters this collection of sets, i.e., for any vector
\smash{$v \in \cb{0, \, 1}^d$}, there exist a policy $\pi_v \in \Pi$ such that $\pi_v(x) = v_j$
for all $x \in \acal_j$. Our proof starts with such a collection of sets \smash{$\cb{\acal_j}_{j = 1}^d$} and
a distribution $\pp$ over $\xx_s$ such that
\begin{equation}
\label{eq:acalvar}
 \EE[\pp]{1\p{\cb{X_i \in \acal_j}} \frac{\sigma^2(X_i)}{e(X_i)(1 - e(X_i))}} = \frac{S_{\pp}}{d} \ \eqfor \ j = 1, \, ..., \, d,
\end{equation}
where $S_{\pp}$ is as defined in \eqref{eq:aabound}. We will establish our result by studying learning
over $\Pi$ with features drawn from this distribution $\pp$.

Now, to lower-bound the minimax risk for policy learning for unknown bounded treatment effect functions
$\tau(\cdot)$, it is sufficient to bound minimax risk over a smaller class of policies $T$, as minimax
risk increases with the complexity of the class $T$. Noting this fact, we restrict our analysis to treatment
functions $T$ such that
$$ \tau(x) = \frac{\sigma^2(x) \, c_j}{e(x)(1 - e(x))} \, \bigg/ \, \EE{\frac{\sigma^2(x) \, 1\p{\cb{X_i \in \acal_j}}}{e(X_i)(1 - e(X_i))}} $$
for all $x \in \acal_j$, where $c_j \in \RR$ is an unknown coefficient for each $j = 1, \, ..., \, d$.
If we knew the values of $c_j$ for $j = 1, \, 2, \,  ..., \, d$, the optimal policy $\pi^* \in \Pi$ would
be treat only those $j$-groups with a positive $c_j$, i.e., $\pi^*(x) = 1\p{\cb{c_j > 0}}$
for all $x \in \acal_j$.

Now, following the argument of \citet{hirano2009asymptotics} (we omit details for brevity),
the minimax policy learner is of the form \smash{$\hpi^*(x) = 1(\{\hc^*_j > 0\})$}
for all $x \in \acal_j$, where \smash{$\hc^*_j$} is an efficient estimator for $c_j$. Moreover, in this
example, we can use \eqref{eq:acalvar} to verify that the semiparametric efficient variance for estimating $c_j$
is $S_{\pp}/d$. Thus, the efficient estimator \smash{$\hc^*_j$} will incorrectly estimate
the sign of $c_j$ with probability tending to \smash{$\Phi(-c_j\sqrt{d/S_{\pp}})$},
where $\Phi(\cdot)$ denotes the standard Gaussian cumulative distribution function.
(Recall that, in our sampling model \eqref{eq:aa}, the signal also decays as $1/\sqrt{n}$.)

By construction, we suffer an expected utility loss of $2\abs{c_j}$ from failing to accurately
estimate the sign of $c_j$. Thus, by the above argument, given fixed values of $c_j$,
the policy learner will suffer an asymptotic regret
$$ \limn \sqrt{n} \EE{R_n} = \sum_{j = 1}^d 2 \abs{c_j} \Phi\p{-\abs{c_j}\sqrt{d/S_{\pp}}}, $$
using an efficient estimator $\hc_j^*$.
Setting $\abs{c_j} = \smash{0.75 \sqrt{S_{\pp}/d}}$, this limit becomes
$$ \limn \sqrt{n} \EE{R_n} = 1.5 \Phi(-0.75) \sqrt{d\, S_{\pp}}, $$
which, noting that $1.5 \Phi(-0.75) \geq 0.33$, concludes the proof.

\end{appendix}

\setlength{\bibsep}{0.2pt plus 0.3ex}
\bibliographystyle{plainnat-abbrev}
\bibliography{references}

\end{document}